\documentclass{elsarticle}

\usepackage{lineno,hyperref}
\usepackage{amsmath,amssymb,amsopn}
\usepackage{url}
\usepackage{graphicx}
\usepackage[caption=false]{subfig}
\usepackage{booktabs}

\usepackage[colorinlistoftodos,prependcaption,textsize=footnotesize]{todonotes}
\reversemarginpar

\usepackage[normalem]{ulem}

\makeatletter
\g@addto@macro\bfseries{\boldmath}
\makeatother

\newtheorem{thm}{Theorem}
\newtheorem{remark}[thm]{Remark}

\numberwithin{equation}{section}

\newcommand{\R}{\mathbb{R}}

\DeclareMathOperator*{\diam}{diam}

\newcommand{\dx}{\,\mathrm{d}x}

\newcommand{\mcK}{\mathcal{K}}

\newcommand{\mcO}{\mathcal{O}}

\newcommand{\hatOmega}{\widehat{\Omega}}
\newcommand{\hatmcK}{\widehat{\mcK}}

\newcommand{\OO}{\mathcal{O}}

\begin{document}

%---------------------------------------------------------------------------
\begin{frontmatter}

\title{\bf MultiMesh Finite Element Methods: \\ Solving PDEs on Multiple Intersecting Meshes}

\author[sintef,simula]{August Johansson\corref{cor1}}
\ead{august.johansson@sintef.no}

\author[simula]{Benjamin Kehlet}
\ead{benjamik@simula.no}

\author[umu]{Mats G. Larson}
\ead{mats.larson@umu.se}

\author[cth]{Anders Logg}
\ead{logg@chalmers.se}

\cortext[cor1]{Corresponding author}

\address[sintef]{SINTEF Digital, Mathematics and Cybernetics, PO Box 124 Blindern, 0314 Oslo, Norway}
\address[simula]{Simula Research Laboratory, PO Box 134, 1325 Lysaker, Norway}
\address[umu]{Department of Mathematics and Mathematical Statistics, Ume{\aa} University, 90187 Ume{\aa}, Sweden}
\address[cth]{Department of Mathematical Sciences, Chalmers University of Technology and University of Gothenburg, 41296 G\"oteborg, Sweden}

\begin{abstract}
  We present a new framework for expressing finite element methods on multiple intersecting meshes: multimesh finite element methods. The framework enables the use of separate meshes to discretize parts of a computational domain that are naturally separate; such as the components of an engine, the  domains of a multiphysics problem, or solid bodies interacting under the influence of forces from surrounding fluids or other physical fields. Such multimesh finite element methods are particularly well suited to problems in which the computational domain undergoes large deformations as a result of the relative motion of the separate components of a multi-body system. In the present paper, we formulate the multimesh finite element method for the Poisson equation. Numerical examples demonstrate the optimal order convergence, the numerical robustness of the formulation and implementation in the face of thin intersections and rounding errors, as well as the applicability of the methodology.

  In the accompanying paper~\cite{mmfem-2}, we analyze the proposed method and prove optimal order convergence and stability.
\end{abstract}


%\begin{AMS}
%  65N30, 65N85,65Y99, 68U20
%\end{AMS}

\end{frontmatter}

%\linenumbers

%---------------------------------------------------------------------------
\section{Introduction}

Finite element methods (FEM) have been successfully applied to almost all areas of science and engineering where the fundamental model is a partial differential equation or system of partial differential equations, and the applicability and generality of finite element methods are undisputed. This is in part due to the generality of the mathematical formulation, allowing stable numerical schemes of arbitrarily high order accuracy to be formulated for all kinds of PDEs. Another contributing factor is the ability of the finite element method to handle computational domains of complex geometry in both 2D and 3D.

However, the wide applicability of FEM relies on the generation of high-quality computational meshes. For simple domains, mesh generation is not an issue but for complex domains, in particular for time-dependent problems or optimization problems with evolving geometries for which a new mesh may need to be generated many times, mesh generation becomes a limiting factor.

\begin{figure}%[htbp]
  \centering{}
  \includegraphics[width=0.4\textwidth]{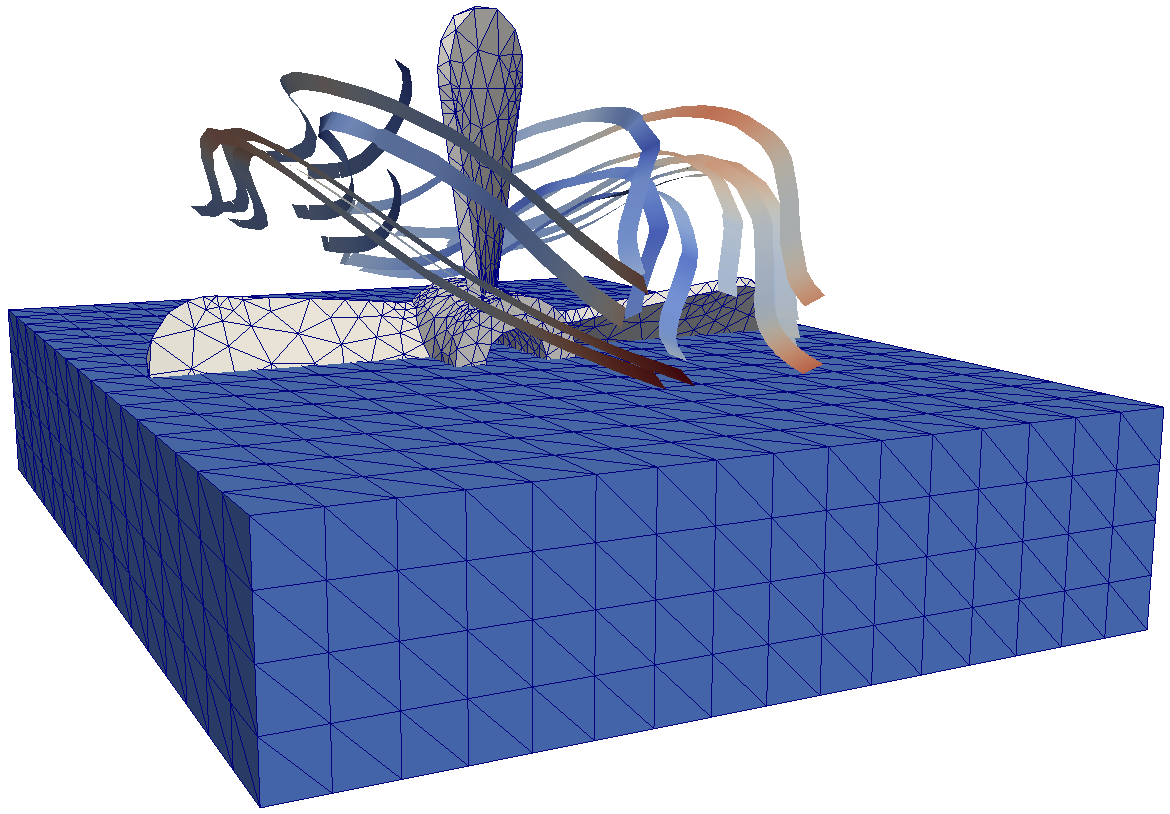} \quad
  \includegraphics[width=0.4\textwidth]{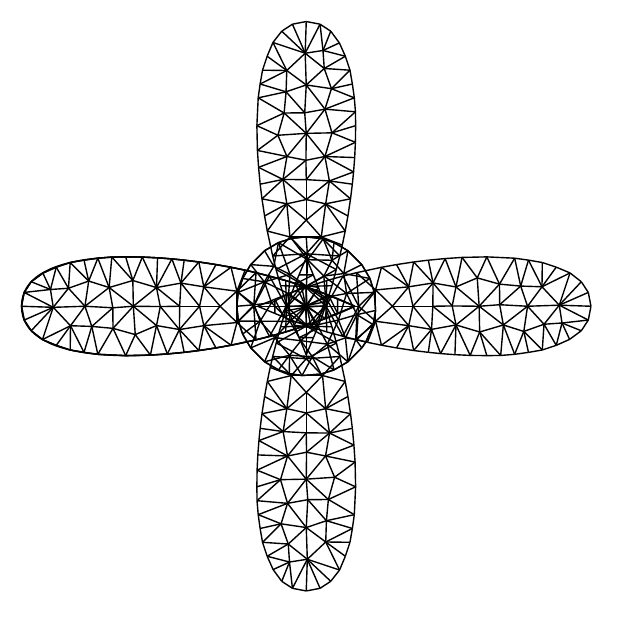} \\[0.5em]
  \includegraphics[width=0.85\textwidth]{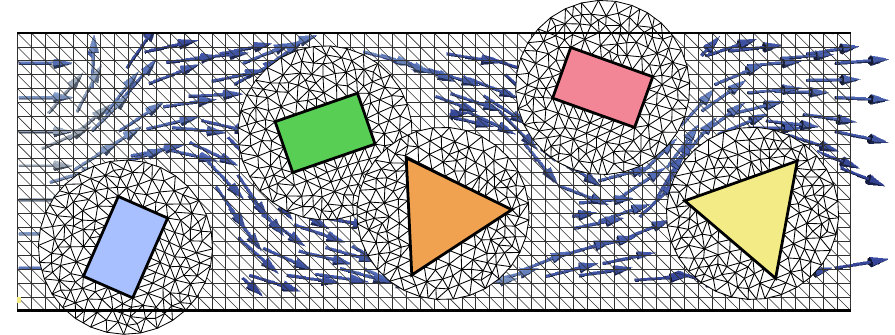}
  \caption{(Top left) The flow around a propeller may be computed by immersing a mesh of the propeller into a fixed background mesh. (Top right) The geometry of a composite object may be discretized by superimposing meshes of each component. (Bottom) The interaction of a set of solid bodies may be simulated using individual meshes that move and intersect freely relative to each other and a fixed background mesh.}
  \label{fig:motivation}
\end{figure}

In the current work, we propose a new way to tackle the mesh generation challenge. By allowing the finite element method to be formulated not on a single mesh but on a multitude of meshes, mesh (re)generation is no longer an issue. For example, a rotating propeller may be approximated by a boundary-fitted mesh that rotates freely on top of a fixed background mesh, the interaction of a set of solid bodies may be simulated using a set of individual meshes that may move and intersect freely during the course of a simulation, and the discretization of a composite object may be obtained by superimposing meshes of its individual components. These cases are illustrated in Figure~\ref{fig:motivation}.

%\pagebreak

In short, multimesh finite element methods extend the finite element method to easily handle computational domains composed of \emph{parts}, where the parts can be either positive or negative (that is, the parts are holes), and where the relative positions and orientations of the parts may change throughout the course of a simulation.

The successful formulation and implementation of multimesh finite element methods rely on three crucial, and challenging, ingredients: (i) consistent finite element formulations that \emph{glue together} the finite element fields expressed on separate meshes on the common interface, (ii) proper stabilization to ensure stability and optimal order convergence in the presence of very thin intersections between elements from separate meshes, and (iii) efficient and robust computational algorithms and implementation that can quickly and robustly compute the intersections between a multitude of arbitrarily intersecting meshes. The latter part is a particular challenge in the face of rounding errors and floating point arithmetic. It should be noted that critical cases will almost surely appear throughout the course of a simulation involving a large number of relative mesh positions, such as in a time-dependent problem.

\subsection{Related work}

The multimesh finite element method is based on Nitsche's method~\cite{Nitsche:1971aa}, which was initially presented as a method for weakly imposing Dirichlet boundary conditions and later extended to form the basis of discontinuous Galerkin methods~\cite{Arnold:1982aa}. Nitsche's method is also the basis for CutFEM, the finite element method on cut meshes. The cut finite element was originally proposed in~\cite{Hansbo:2002aa,Hansbo:2003aa} and has since been applied to a range of problems~\cite{Burman:2007aa,Burman:2007ab,Becker:2009aa,Massing:2015aa}. For an overview, see~\cite{Burman:2015aa} and~\cite{Bordas-et-al-2017}.

The problem of handling non-matching and multiple meshes is discussed in the literature under many different names, such as fictitious domain, unfitted methods, overset grids and overlapping meshes. Numerous numerical approaches also exists, for example using Lagrange multipliers~\cite{BURMAN20102680}, finite cell methods~\cite{Schillinger2015,DUSTER20083768,Parvizian2007}, immersed interface methods~\cite{LI1998253}, the s-version of the finite element method \cite{FISH1992539,FISH1994135} and XFEM~\cite{NME:NME2914,NME:NME386,LEHRENFELD2016716} to name a few. Discontinuous Galerkin methods have also successfully been used~\cite{NME:NME2631,Johansson:2013:HOD:2716605.2716805,SAYE2017647,SAYE2017683}, as well as variants using local enrichments~\cite{BOUCLIER20161,Zander2015,SCHILLINGER2012116}. There are also substantial and important works in domain decomposition methods see for example~\cite{RANK1992299,BeckerHansboStenberg} or \cite{APPELO20126012,Henshaw:2008:PCT:1387360.1387531} and the references therein. There are also techniques useful in particular cases, for example in the case of a sliding mesh on top of a fixed background mesh \cite{behr-tezduyar,bazilevs_korobenko}.

In cut finite element methods, the finite element solution is represented as a standard (continuous piecewise polynomial) finite element function, either on a fixed background mesh or on a pair of intersecting meshes. Boundary conditions are then imposed weakly either on an intersecting interface, defined for example by a level set function or a discrete surface, or on the interface between two meshes to ensure near continuity of the global finite element field. These boundary conditions are imposed using Nitsche's method and the method requires the addition of suitable stabilization terms, usually in the form of penalizing jumps in function values or derivatives across interfaces, to ensure stability and optimal order convergence.

The multimesh finite element method is a generalization of CutFEM to general collections of overlapping meshes, where three or more meshes may overlap simultaneously. This imposes new requirements for stabilization and, notably, brings new challenges for robust algorithms and implementation in the face of rounding errors and floating point arithmetic, when tens or hundreds of meshes must be correctly intersected at once, and quadrature rules be computed on the complex shapes resulting from mesh intersections.

The current work is closely related to the work presented in~\cite{Massing:2014aa,Johansson:2015aa} where two overlapping meshes were used to discretize the Stokes problem. The algorithms are similar to those presented in~\cite{Massing:2013ab} in that axis-aligned bounding boxes (AABB trees) are used for efficient mesh intersection, but different in that we must now consider the intersection of $N+1$ meshes, not just a pair of meshes, and that a new algorithm has been used to compute quadrature rules on cut cells and interfaces. For a detailed exposition of these algorithms and software, we refer to the related work~\cite{Johansson:2017ab}.

\subsection{Outline}

The remainder of this paper is organized as follows. The multimesh finite element method for the Poisson problem is presented in Section~\ref{sec:method}. The algorithms and implementation of the multimesh finite element method are then briefly discussed in Section~\ref{sec:algorithms}. Numerical results are presented in Section~\ref{sec:numres}. Finally, conclusions are presented in Section~\ref{sec:conclusions}.

%---------------------------------------------------------------------------
\section{The MultiMesh Finite Element Method}
\label{sec:method}

The starting point for the multimesh finite element method is a sequence of standard finite element spaces on $1 + N$ meshes, consisting of a single background mesh covering the computational domain and $N$ overlapping/intersecting meshes arbitrarily positioned on top of the background mesh.
The finite element solution is then expressed as a composite function in the multimesh finite element space and is glued together using Nitsche's method with appropriate stabilization.

Before we express the method in detail, we here first introduce the concepts and notation for domains, interfaces, meshes, overlaps and function spaces.

\subsection{Domains}
\label{sec:domainsetup}

\begin{figure}
  \centering
  \subfloat[]{\label{fig:three_domains_a}\includegraphics[height=0.25\linewidth]{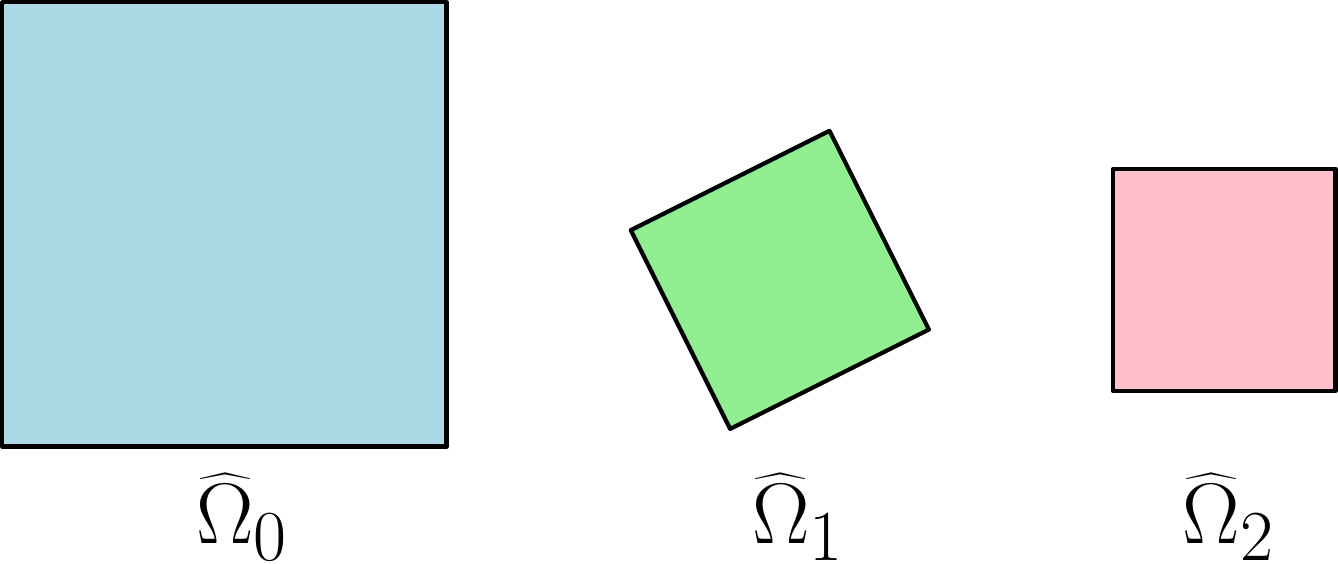}}\qquad\qquad\qquad
  \subfloat[]{\label{fig:three_domains_b}\includegraphics[height=0.25\linewidth]{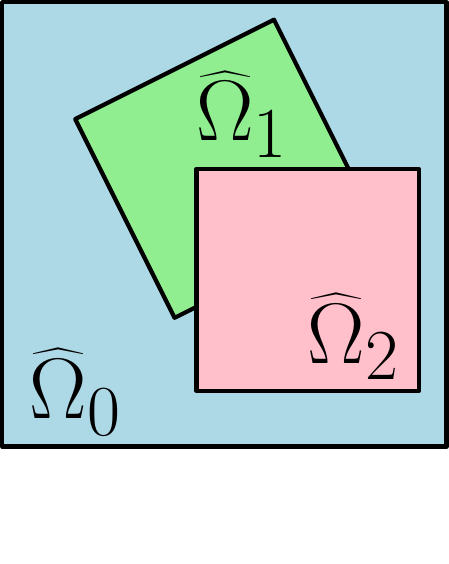}}\\
  \caption{(a) Three polygonal predomains. (b) The predomains are placed on top of each other in an ordering such that
    $\hatOmega_0$ is placed lowest, $\hatOmega_1$ is in the middle and $\hatOmega_2$ is on top. }
  \label{fig:three_domains}
\end{figure}

\begin{figure}
  \centering
  \includegraphics[height=0.2\linewidth]{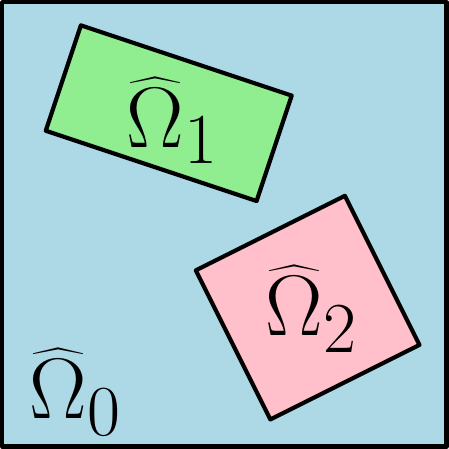}
  \caption{Three predomains in a configuration where the order of $\hatOmega_1$ and $\hatOmega_2$ is not unique. Swapping the order of $\hatOmega_1$ and $\hatOmega_2$ will lead to the same multimesh.}
  \label{fig:non_unique}
\end{figure}

Let $\Omega = \hatOmega_0 \subset \R^d$, $d = 2,3$, be a domain with polygonal boundary (the background domain) and let
\begin{align}
  \hatOmega_i \subset \hatOmega_0, \qquad i=1,\ldots, N,
\end{align}
be subdomains of $\hatOmega_0$ with polygonal boundaries. We impose an ordering of the domains by saying that $\hatOmega_j$ is on top of $\hatOmega_i$ if $j > i$. For an illustration, see Figure~\ref{fig:three_domains}. Note that the ordering is not unique as illustrated in Figure~\ref{fig:non_unique}. We refer to the domains $\hatOmega_0,\ldots,\hatOmega_N$ as the \emph{predomains}. The predomains may or may not intersect each other, but they will intersect the background domain.

\begin{figure}
  \begin{center}
    \includegraphics[height=0.25\linewidth]{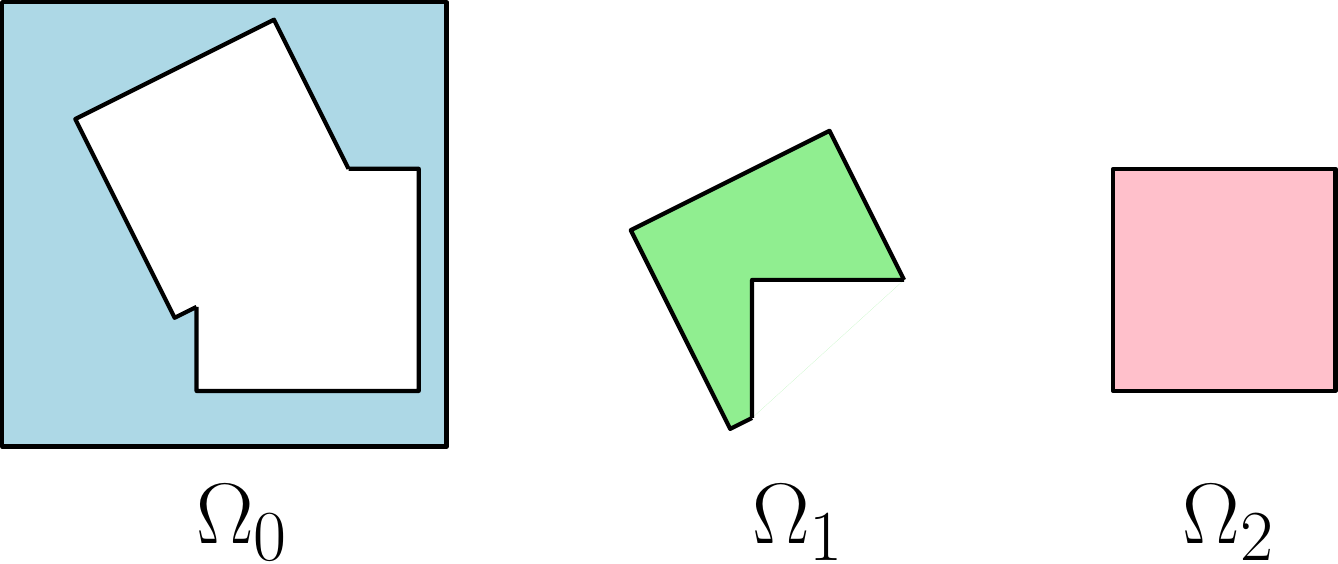}
    \caption{Partition of $\Omega = \Omega_0 \cup \Omega_1 \cup \Omega_2$. Note that $\Omega_2 = \hatOmega_2$.}
    \label{fig:three_domains_partition}
  \end{center}
\end{figure}

We now define a partition $\{ \Omega_i \}_{i=0}^N$ of $\Omega$ by letting
\begin{align}
  \label{eq:Omega_i}
  \Omega_i = \hatOmega_i \setminus \bigcup_{j=i+1}^{N} \hatOmega_j,  \qquad i=0, \ldots, N.
\end{align}
In other words, \emph{the domain $\Omega_i$ is the visible part of the predomain $\hatOmega_i$}; see Figure~\ref{fig:three_domains_partition}. Note that by definition, we have have $\Omega_N = \hatOmega_N$. Also note that some domains may be completely hidden behind domains with a higher index, and thus be empty.

\begin{remark}
  \label{rem:boundary-overlap}
By assumption, the domains $\Omega_1, \ldots, \Omega_N$ are not allowed to intersect the boundary of $\Omega$. This assumption may be too restrictive for some applications, such as for two of the domains shown in the bottom of Figure~\ref{fig:motivation}, but we make this assumption for simplicity in the current manuscript. The multimesh finite element method, and indeed also our implementation, support computational domains that extend beyond the boundary of the background mesh.
\end{remark}

\subsection{Interfaces}

\begin{figure}
  \centering
  \subfloat[]{\label{fig:two_interfaces_a}\includegraphics[width=0.25\linewidth]{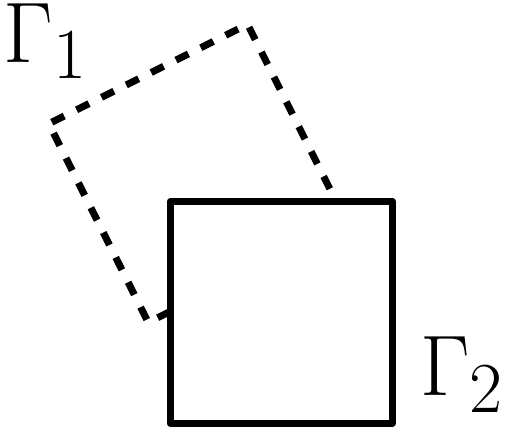}}\qquad\qquad\qquad
  \subfloat[]{\label{fig:two_interfaces_b}\includegraphics[width=0.25\linewidth]{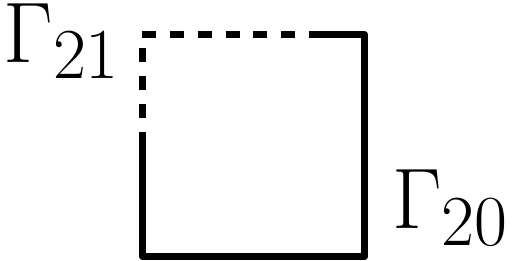}}
  \caption{(a) The two interfaces of the domains in Figure~\ref{fig:three_domains}: $\Gamma_1 = \partial \hatOmega_1 \setminus \hatOmega_2$ (dashed line) and  $\Gamma_2 = \partial \hatOmega_2$ (filled line). Note that $\Gamma_1$ is not a closed curve. (b) By~\eqref{eq:gamma_partition}, $\Gamma_1 = \Gamma_{10}$ and $\Gamma_2 = \Gamma_{20} \cup \Gamma_{21}$.}
\end{figure}

Given a sequence of predomains $\{\hatOmega_i \}_{i=0}^N$, let
\begin{align}
  \Gamma_i = \partial \hatOmega_i \setminus \bigcup_{j=i+1}^N \hatOmega_j,  \qquad i=1, \ldots, N-1.
\end{align}
In other words, the \emph{interface} $\Gamma_i$ is the visible part of the boundary of $\hatOmega_i$; see Figure~\ref{fig:two_interfaces_a}. Note that by definition, we have $\Gamma_N = \partial \hatOmega_N$. Also note that $\Gamma_i$ does not need to be a closed curve. See also Remark~\ref{rem:boundary-overlap}.

We may further partition each interface $\Gamma_i$ into a set of disjoint \emph{interfaces},
\begin{align}
  \label{eq:gamma_partition}
  \Gamma_i = \bigcup_{j=0}^{i-1} \Gamma_{ij}
           = \bigcup_{j=0}^{i-1} \Gamma_i \cap \Omega_j,
  \qquad i=1, \ldots, N.
\end{align}
In other words, the interface $\Gamma_{ij} = \Gamma_i \cap \Omega_j$ is the subset of the visible boundary of the predomain $\hatOmega_i$ that intersects with the domain $\Omega_j$, where $j < i$; see Figure~\ref{fig:two_interfaces_b}. Note that some $\Gamma_{ij}$ may be empty.

\subsection{Meshes}

\begin{figure}
  \centering
  \subfloat[]{\label{fig:three_meshes_a}\includegraphics[height=0.25\linewidth]{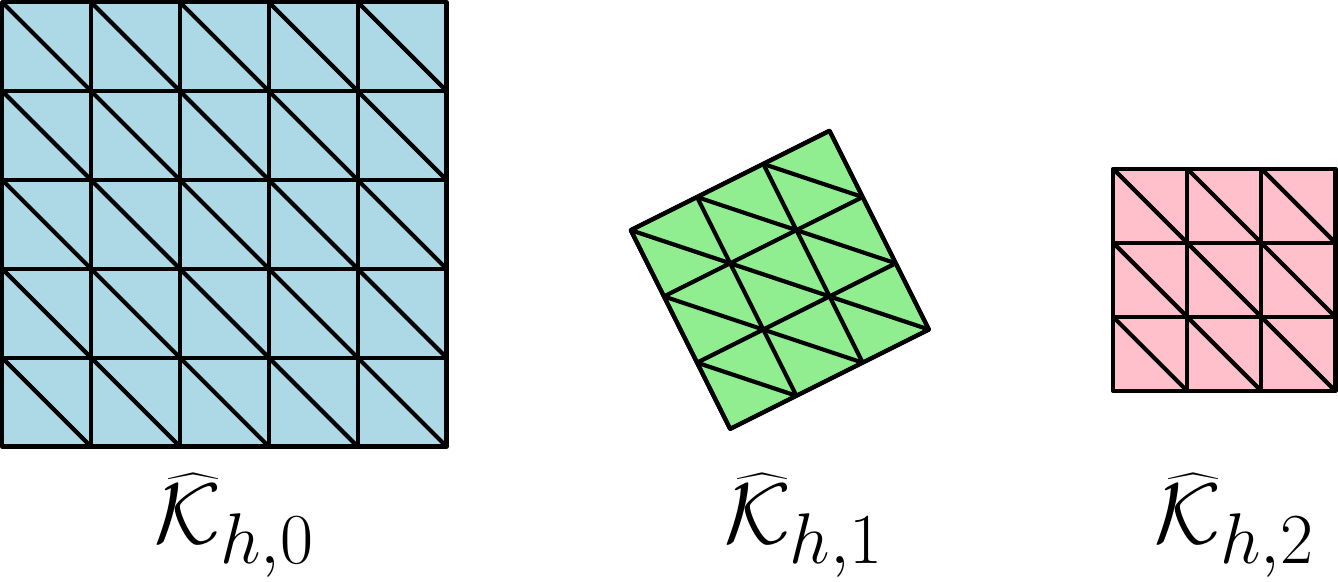}}\qquad\qquad\qquad
  \subfloat[]{\label{fig:three_meshes_b}\includegraphics[height=0.25\linewidth]{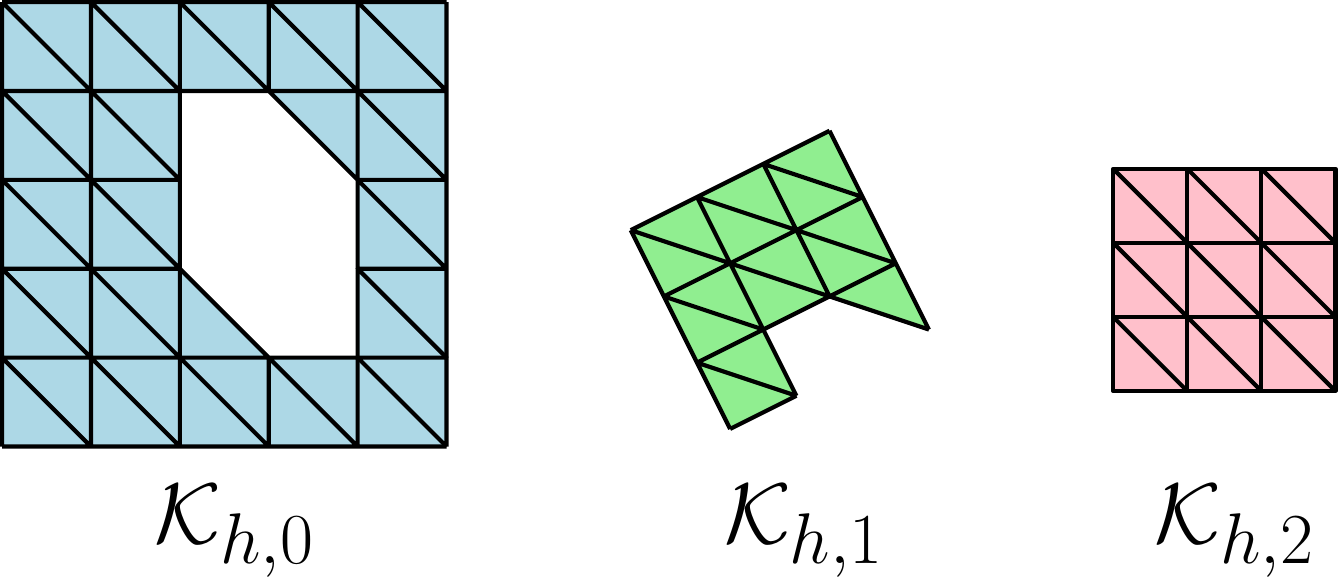}}
  \caption{(a) The three premeshes. (b) The corresponding active meshes (cf.\ Figure~\ref{fig:three_domains}).}
\end{figure}

Recall that for a quasi-uniform triangular mesh $\mcK_h$ there exist constants $c_0, c_1$ such that for all $K\in \mcK_h$ (see e.g.~\cite{BreSco08}),
\begin{align}
  \max_{K\in \mcK_h} \frac{\diam(K)^d}{|K|} &\leq c_0, \qquad\qquad  \frac{\max_{K\in \mcK_h} |K|}{\min_{K\in \mcK_h} |K|} \leq c_1.
\end{align}
where $|K|$ denotes the Lebesgue measure of $K$.
Let now $\hatmcK_{h,i}$ be a quasi-uniform \emph{premesh} on $\hatOmega_i$ with mesh parameter
\begin{align}
  h_i = \max_{K\in \hatmcK_{h,i}} \diam(K), \qquad i=0,\ldots,N.
\end{align}
In other words, the premesh $\hatmcK_{h,i}$ is a mesh of the predomain $\hatOmega_i$ created by a standard mesh generator; see Figure~\ref{fig:three_meshes_a}. We define the global mesh size $h$ by $h = \max_i h_i$.

\begin{figure}
  \centering
  \subfloat[]{\label{fig:overlap}\includegraphics[height=0.25\linewidth]{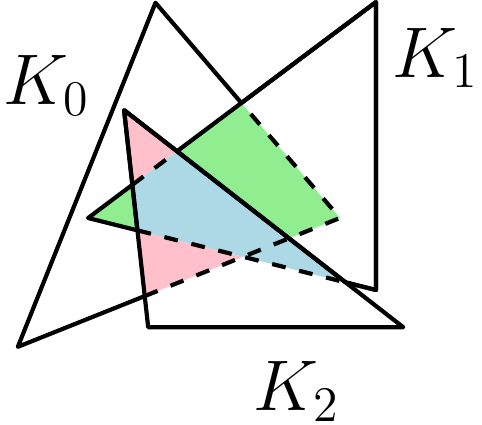}}\qquad\qquad\qquad
  \subfloat[]{\label{fig:multimesh}\includegraphics[height=0.25\linewidth]{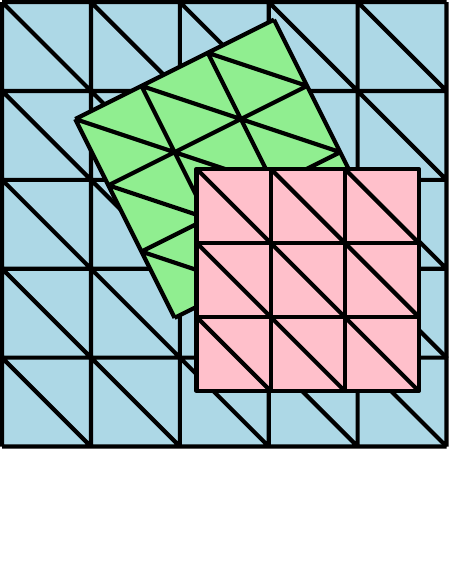}}
  \caption{(a) Given three ordered triangles $K_0$, $K_1$ and $K_2$, the overlaps are $\mcO_{01}$ in green, $\mcO_{02}$ in red and $\mcO_{12}$ in blue. (b) The multimesh of the domains in Figure~\ref{fig:three_domains_b} consists of the active meshes in Figure~\ref{fig:three_meshes_b}.}
\end{figure}

We further define the \emph{active meshes} by
\begin{align}  \label{eq:active_mesh}
  \mcK_{h,i} = \{ K \in \hatmcK_{h,i} : K \cap \Omega_i \neq \emptyset \},
  \qquad i=0,\ldots,N.
\end{align}
In other words, the active mesh $\mcK_{h,i}$ consists of all (fully or partly) visible elements of the premesh $\hatmcK_{h,i}$; see Figure~\ref{fig:three_meshes_b}. Placed in the given ordering, the active meshes form the multimesh; see Figure~\ref{fig:multimesh}. Note that the active mesh consists of regular (uncut) elements as well as elements which may be partly hidden behind overlapping meshes (the cut elements).

We also define the \emph{active domains} by
\begin{align}
  \Omega_{h,i} = \bigcup_{K\in\mcK_{h,i}} K,
  \qquad i=0,\ldots,N.
\end{align}
In other words, the active domain $\Omega_{h,i}$ is the domain defined by the active mesh $\mcK_{h,i}$. Note that $\Omega_{h,i}$ typically extends beyond the corresponding domain $\Omega_i$.

\subsection{Overlaps}

Given a sequence of active domains $\{\Omega_{h,i}\}_{i=0}^N$, let
\begin{align}
  \OO_i = \Omega_{h,i} \setminus \Omega_i,
  \qquad i=0,\ldots,N-1.
\end{align}
In other words, the \emph{overlap} $\OO_i$ is the hidden part of the active domain $\Omega_{h,i}$. We may further partition each overlap into a set of disjoint \emph{overlaps},
\begin{align}
  \OO_i = \bigcup_{j=i+1}^N \OO_{ij}
        = \bigcup_{j=i+1}^N \OO_i \cap \Omega_j,
  \qquad i=0, \ldots, N-1.
\end{align}
In other words, the overlap $\OO_{ij} = \OO_i \cap \Omega_j = \Omega_{h,i} \cap \Omega_j$ is the hidden part of the elements of the mesh $\mcK_{h,i}$ that are hidden below the domain $\Omega_j$. Since some $\OO_{ij}$ may be empty, we will make use of the indicator function $\delta_{ij}$  defined by
\begin{align}
  \label{eq:indicatorfunction}
  \delta_{ij} =
  \begin{cases}
    1, \quad \OO_{ij} \neq \emptyset,
    \\
    0, \quad \text{otherwise},
  \end{cases}
\end{align}
and let $N_\OO$ denote the maximum number of non-empty overlaps
\begin{align}
  \label{eq:max-overlaps}
  N_\OO = \max_{1 \leq j \leq N} \sum_{i=0}^{N-1} \delta_{ij}.
\end{align}
Note that $N_\OO \leq N$. For example, in Figure~\ref{fig:three_domains}, $N_\OO = N = 3$ since all three predomains intersect. In Figure~\ref{fig:non_unique}, $N = 3$ and $N_\OO = 2$ since there are only two intersecting domains.

\subsection{Function spaces}
\label{sec:fespaces}

Given a sequence of active meshes $\{\mcK_{h,i}\}_{i=0}^N$, let
$V_{h,i}$ be a continuous piecewise polynomial finite element space on the mesh $\mcK_{h,i}$, the active elements of domain $i$. Note that even if only a part of an element in $\mcK_{h,i}$ is visible, the domain extends to the whole elements.

We now define the multimesh finite element space as the direct sum of the individual finite element spaces $\{V_{h,i}\}_{i=0}^N$:
\begin{align}
  V_h = \bigoplus_{i=0}^N V_{h,i}.
\end{align}
An element $v$ of $V_h$ is a tuple $(v_0, \ldots, v_N)$ and the inclusion $V_h \hookrightarrow L^2(\Omega)$ is defined by the evaluation
\begin{align}
  v(x) = v_i(x),
  \quad x\in\Omega_i.
\end{align}
In other words, $v(x)$ is evaluated by evaluating $v_i(x)$ for the highest ranked (topmost) mesh containing the point $x$.

\subsection{Finite element method}
\label{sec:formulation}

We may now formulate the multimesh finite element method for the Poisson problem
\begin{align}
  \label{eq:poisson}
  -\Delta u = f \quad \text{in } \Omega,
  \qquad u = 0 \quad \text{on } \partial \Omega,
\end{align}
where $\Omega \subset \R^d$ is a polygonal domain.

Given the construction of the finite element spaces in Section~\ref{sec:fespaces}, the multimesh finite element method for~\eqref{eq:poisson} is to find $u_h \in V_h$ such that
\begin{align}
  \label{eq:method}
  A_{h}(u_h, v)  &= l_h(v) \qquad \forall v \in V_h,
\end{align}
where
\begin{align}\label{eq:Ah}
A_h(v,w) & =  a_{h}(v, w) + s_h(v, w),
\\ \label{eq:ah}
  a_{h}(v, w) &=
   \sum_{i=0}^N (\nabla v_i,\nabla w_i)_{\Omega_i}
   \\ \nonumber
   &\qquad
   - \sum_{i=1}^N \sum_{j=0}^{i-1} (\langle n_i\cdot \nabla v \rangle, [w])_{\Gamma_{ij}}
      + ([v], \langle n_i \cdot \nabla w \rangle)_{\Gamma_{ij}}
\\ \nonumber
&\qquad
 + \sum_{i=1}^N \sum_{j=0}^{i-1} \beta_0 ( h_i^{-1} + h_j^{-1} ) ([v], [w])_{\Gamma_{ij}},
  \\ \label{eq:sh}
  s_h (v, w) &= \sum_{i=0}^{N-1} \sum_{j=i+1}^N \beta_1 ([ \nabla v ], [\nabla w])_{\OO_{ij}},
  \\ \label{eq:lh}
  l_h(v) &= \sum_{i=0}^N (f,v_i)_{\Omega_i}.
\end{align}
Here, $\beta_0 > 0$ and $\beta_1 > 0$ are a stabilization parameters that must be sufficiently large to ensure that the bilinear form $A_h$ is coercive; see~\cite{mmfem-2} for an analysis.

Note that the sum of interface terms extends over the range $0 \leq j < i \leq N$, since the interface $\Gamma_i$ is partitioned into interfaces $\Gamma_{ij}$ relative to underlying meshes, whereas the sum of overlap terms extends over the range $0 \leq i < j \leq N$, since the overlap $\OO_i$ is partitioned into overlaps $\OO_{ij}$ relative to overlapping meshes.

We note that the bilinear form $a_h$ contains the standard symmetric Nitsche formulation of~\eqref{eq:poisson} with jumps and fluxes given pairwise over $\Gamma_{ij}$. The contribution from the bilinear form $s_h$ is a stabilization of the jump in the gradients over $\mcK_{h,i} \cap \mcK_{h,j}$. The effect of this term is analyzed in~\cite{mmfem-2}. In short, it ensures that the bilinear form $A_h$ is coercive, which guarantees a unique solution of the problem by the Lax-Milgram theorem~\cite{BreSco08}. The term also ensures the conditioning of the linear system, which otherwise scales as $\eta^{-(2p+1-2/d)}$~\cite{Prenter-et-al}, where $\eta$ is the volume fraction of the smallest cell and $p$ is the polynomial order. By including $s_h$, the condition number scales independently of the volume fraction and only with respect to the mesh size. This is also seen in \cite{Burman:2015aa,Johansson:2017ab} and similar stabilized methods.

The jump terms on $\Gamma_{ij}$ and $\OO_{ij}$ are given by $[v] = v_i - v_j$, where $v_i$ and $v_j$ are the finite element solutions on represented on the active meshes $\mcK_{h,i}$ and $\mcK_{h,j}$.
On $\Gamma_{ij}$ we also define the average normal flux by
\begin{equation}\label{eq:average}
\langle n_i\cdot \nabla v \rangle = (n_i \cdot \nabla v_{i} + n_{i} \cdot \nabla v_{j})/2.
\end{equation}
In the average any convex combination is valid~\cite{Hansbo:2003aa}. Note that the finite element method approximates $[u_h] = 0$ and $[n_i \cdot \nabla u_h] = 0$ on all interfaces.

%---------------------------------------------------------------------------
\section{Algorithms and Implementation}
\label{sec:algorithms}

The multimesh finite element method has been implemented as part of FEniCS~\cite{Logg:2012aa,Alnaes:2015aa} and the numerical
examples in Section~\ref{sec:numres} were performed using (a development version of) the 2017.2 release. In this section, we comment briefly on the algorithms and the implementation and refer to the related work~\cite{Johansson:2017ab} for a more detailed exposition.

The evaluation of the integrals appearing in~\eqref{eq:ah},~\eqref{eq:sh} and~\eqref{eq:lh} impose several challenges and rely heavily on the following functionality.
\begin{itemize}
\item \textbf{Collision detection.}
For each mesh, compute which elements in the mesh collide with which elements in all other meshes. This step involves computing $N + (N-1) + \ldots + 1 \approx N^2/2$ mesh--mesh collisions.
\item \textbf{Intersection construction.}
For each mesh and each colliding element, construct the intersection of the element with the colliding mesh. This step may be reduced to repeated triangulation of simplex--simplex intersections.
\item \textbf{Quadrature rule construction.}
For each cut element (not necessarily convex) and each interface, construct a quadrature rule for integration over the cut element or interface.
\end{itemize}

Efficient construction of quadrature rules is an essential component of a method handling non-matching and multiple meshes. Previous approaches include subtriangulation~\cite{Massing:2013ab}, moment-fitting~\cite{Sudhakar:2013aa,Muller:2013aa}, recursive bisectioning~\cite{VERHOOSEL2015138,Schillinger2015}, tensor-product based~\cite{Saye-high-order}, to name a few. One important fact to note in our setting is that curved interfaces can in general be avoided, since the interfaces are completely artificial.

\subsection{Collision detection}

Given a pair of (pre)meshes $\hatmcK_{h,i}$ and $\hatmcK_{h,j}$, we compute the collisions between all elements in $\hatmcK_{h,i}$ with all elements in $\hatmcK_{h,j}$. The collisions are computed by first constructing and then colliding axis-aligned bounding box trees (AABB trees) for the two meshes~\cite{RTR3,Ericson:2004aa} to find candidate element intersections. The candidate elements are then checked for collision using robust geometric predicates of adaptive precision~\cite{Shewchuk:1997aa}. The mesh--mesh collision is carried out for each (unordered) pair of premeshes among the $1 + N$ meshes.

\subsection{Intersection construction}

Given a pair of colliding elements $K_i$ and $K_j$ belonging to two different meshes, we compute a subtriangulation of the intersection $K_i \cap K_j$. Since $K_i \cap K_j$ is convex, the subtriangulation may be computed by a simple Graham scan~\cite{Graham72}. When an element collides with more than one element from the same mesh; that is, if $K_i \in \hatmcK_{h,i}$ collides with $K_j,K_j'\in\hatmcK_{h,j}$, a subtriangulation is created for each pairwise intersection $K_i \cap K_j$ and $K_i \cap K_j'$, etc.

When an element collides with more than one element from multiple meshes; that is, if $K_i \in \hatmcK_{h,i}$ collides with $K_j \in \hatmcK_{h,j}$ and $K_k \in \hatmcK_{h,k}$, a subtriangulation is first created for the intersection $K_i \cap K_j$ (assuming that $j < k$), and then each element in this subtriangulation is further intersected with $K_k$ and new subtriangulations constructed. This procedure guarantees that the only intersections that need to be constructed (triangulated) are pairwise intersections of elements (simplices). This procedure is also the basis for the construction of quadrature rules as described below.

\subsection{Quadrature rule construction}

Several techniques have been developed for integration over polygonal and polyhedral domains; see, e.g.,~\cite{Mousavi:2010aa,Muller:2013aa,Sudhakar:2013aa}. However, the construction and representation of general polygonal and polyhedral domains imposes several algorithmic and implementational challenges. Instead, we propose a new approach which, to the authors knowledge, has not been presented before in the literature. We here give a short overview of the basic (and simple) idea and present the details of the algorithm in~\cite{Johansson:2017ab}.

\begin{figure}
  \centering
  \subfloat[]{\label{fig:cut_cutting_example}\includegraphics[height=0.25\linewidth]{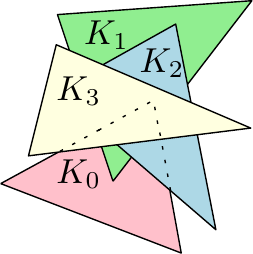}}\qquad\qquad\qquad
  \subfloat[]{\label{fig:integral_types}\includegraphics[height=0.25\linewidth]{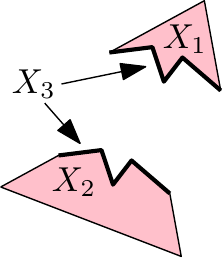}}
  \caption{(a) An element $K_0$ overlapped by the elements $K_1$, $K_2$ and $K_3$.
    (b) The configuration in (a) gives rise to three integration domains with respect to $K_0$, namely $X_1$, $X_2$ and $X_3$.}
  \label{fig:Xs}
\end{figure}

We first note that the multimesh method relies on computing integrals over three kinds of domains; these domains are marked as $X_1$, $X_2$ and $X_3$ in Figure~\ref{fig:Xs}. The domain $X_1$ is the hidden part of a cut cell, the domain $X_2$ is the visible part of the cut cell, and $X_3$ is an interface. In 2D the main challenge is the integration over the polygonal domains $X_1$ and $X_2$ and in 3D the challenge is the corresponding integration over polyhedral domains. For the current discussion we thus focus on the integration over these domains and omit the discussion on the integration over $X_3$ (which is simpler).

Consider first the integration over $X_1$ in Figure~\ref{fig:Xs} and note that the domain is the result of repeatedly intersecting a set of triangles. The result of an intersection between two triangles is always a convex polygon. Since a convex polygon can be trivially triangulated (partitioned into triangles), we can easily represent the intersection of two triangles as a union of triangles. We thus first intersect $K_0$ and $K_1$ in Figure~\ref{fig:Xs} and get a set of (two) triangles. We then continue to intersect these triangles with $K_2$ and $K_3$ and the result is a new set of triangles. The quadrature rule over $X_1$ is then easily constructed by summing any standard quadrature rules over the set of triangles.

\begin{figure}
  \begin{center}
    \includegraphics[width=0.45\linewidth]{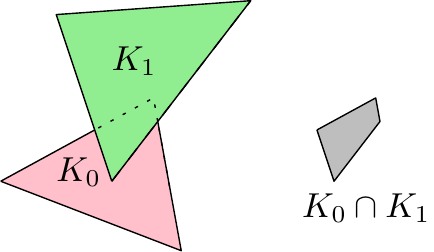}
    \caption{Two triangles $K_0$ and $K_1$ and the convex intersection $K_0 \cap K_1$.}
    \label{fig:Venn-2-tris}
  \end{center}
\end{figure}

Consider next the more challenging integration over $X_2$. If we have only two meshes, we only need to intersect two triangles ($K_0$ and $K_1$). In this case, $X_1$ is given by $K_0\setminus K_1$; see Figure~\ref{fig:Venn-2-tris}.
The key is now to \emph{not} represent $X_1$ but to instead use the well known inclusion-exclusion principle from combinatorics which states that

\begin{equation}
  |K_0 \cup K_1| = |K_0| + |K_1| - |K_0 \cap K_1|.
\end{equation}
In other words, a quadrature rule for the union $K_0 \cup K_1$ can be constructed by summing the quadrature rules for $K_0$ and $K_1$ and subtracting the quadrature rule for the intersection $K_0 \cap K_1$. It follows that a quadrature rule for the difference $K_0 \setminus K_1$ is given by
\begin{align}
  |K_0 \setminus K_1|
  &= (K_0 \cup K_1) \setminus K_1 \\
  &= |K_0| + |K_1| - |K_0 \cap K_1| - |K_1| \\
  &= |K_0| - |K_0 \cap K_1|.
\end{align}
Since now both $K_0$ and $K_0 \cap K_1$ can be easily represented as a set of triangles, we can easily compute a quadrature for $K_0 \setminus K_1$. Since the inclusion-exclusion principle generalizes to the union of arbitrarily many domains, for example
\begin{align}
  |K_0 \cup K_1 \cup K_2|
  &= |K_0| + |K_1| + |K_2|
  \nonumber \\
  &\quad - |K_0 \cap K_1| - |K_0 \cap K_2| - |K_1 \cap K_2| + |K_0 \cap K_1 \cap K_2|,
\end{align}
a quadrature rule for the domain $X_2$ can be easily constructed by adding and subtracting the quadrature rules for a set of triangles. The result is a set of quadrature points and weights for $X_1$, where some of the weights are negative.

%---------------------------------------------------------------------------
\section{Numerical Results}
\label{sec:numres}

In this section, we present a series of numerical results to
demonstrate the accuracy and stability of the multimesh finite element
formulation presented in Section~\ref{sec:formulation}, as well as the
robustness and generality of the implementation.

\subsection{Optimal convergence}
\label{sec:convergence}

We solve the Poisson problem on a discretization of $1 + N$ (pre)domains $\{\hatOmega_i \}_{i=0}^N$ with $\Omega_0 = \hatOmega_0 = [0,1]^2$. The $N$ overlapping domains are placed at (uniformly) random positions and scaled to a (uniformly) random size inside $\Omega_0$ such that they do not intersect with the boundary of $\Omega_0$. We take $\beta_0=6p^2$, where $p$ is the polynomial order, and $\beta_1=10$ and investigate the convergence for $N \in \{0,1,2,4,8,16,32 \}$. As a test problem, we consider the analytical solution
\begin{align}
  \label{eq:exactsolution}
  u(x,y) = \sin(\pi x) \sin(\pi y),
\end{align}
corresponding to the right-hand side $f(x, y) = 2\pi^2 \sin(\pi x) \sin(\pi y)$.

\begin{figure}
  \begin{center}
    \includegraphics[width=0.320\linewidth]{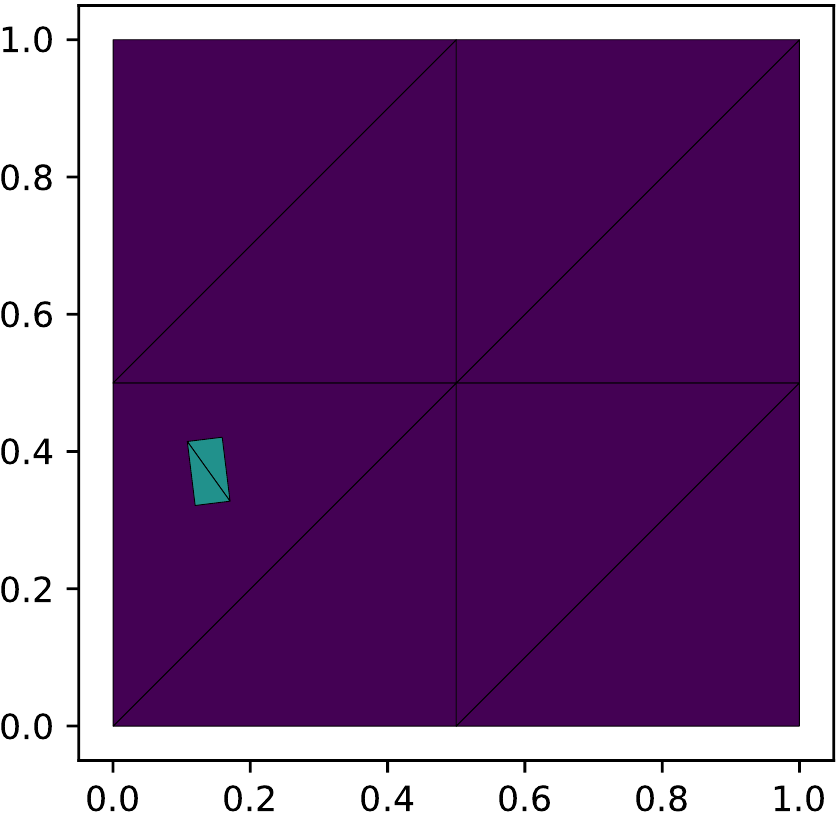}\ \
    \includegraphics[width=0.320\linewidth]{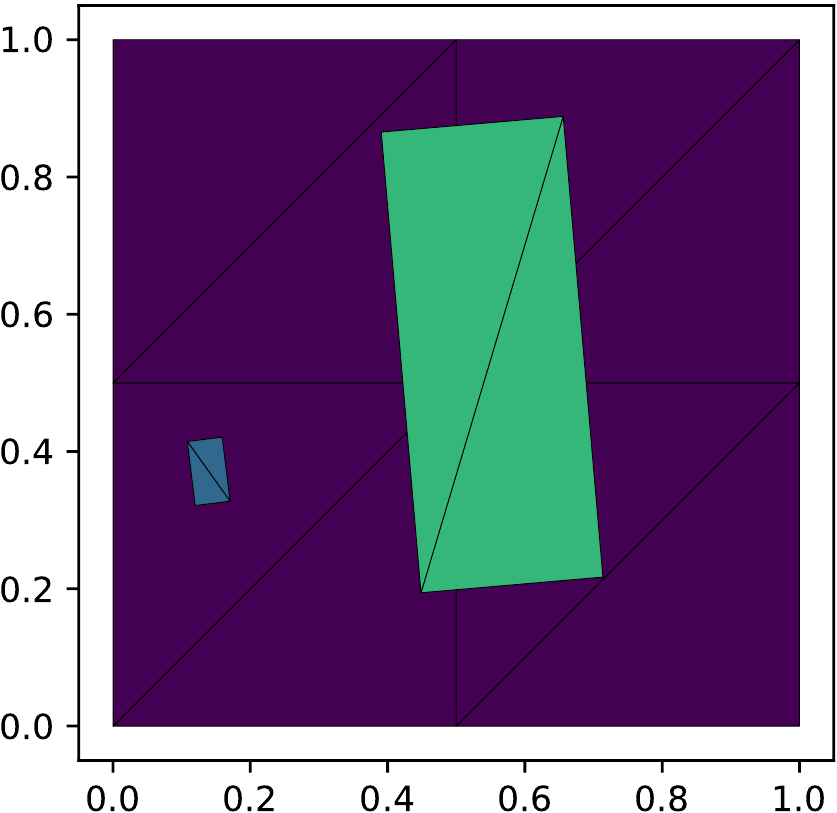}\ \
    \includegraphics[width=0.320\linewidth]{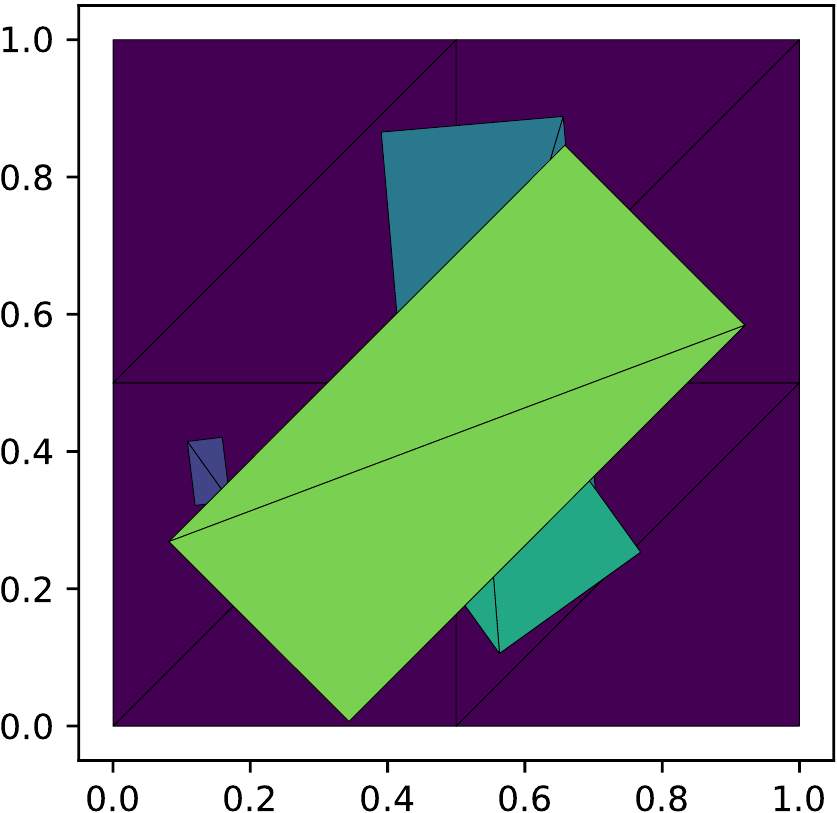}\\ \ \\
    \includegraphics[width=0.320\linewidth]{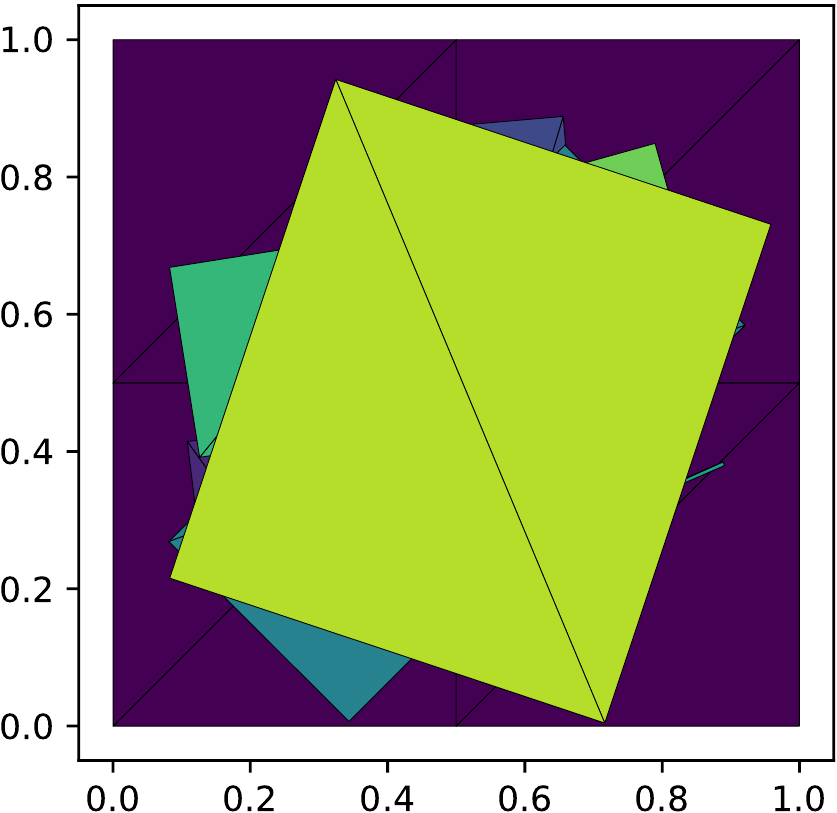}\ \
    \includegraphics[width=0.320\linewidth]{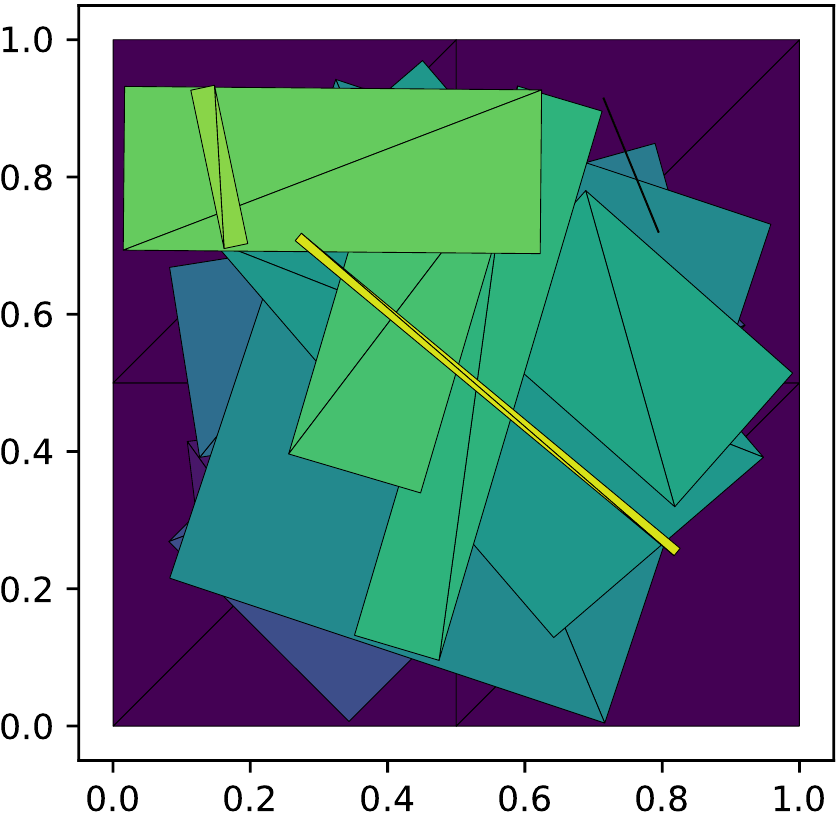}\ \
    \includegraphics[width=0.320\linewidth]{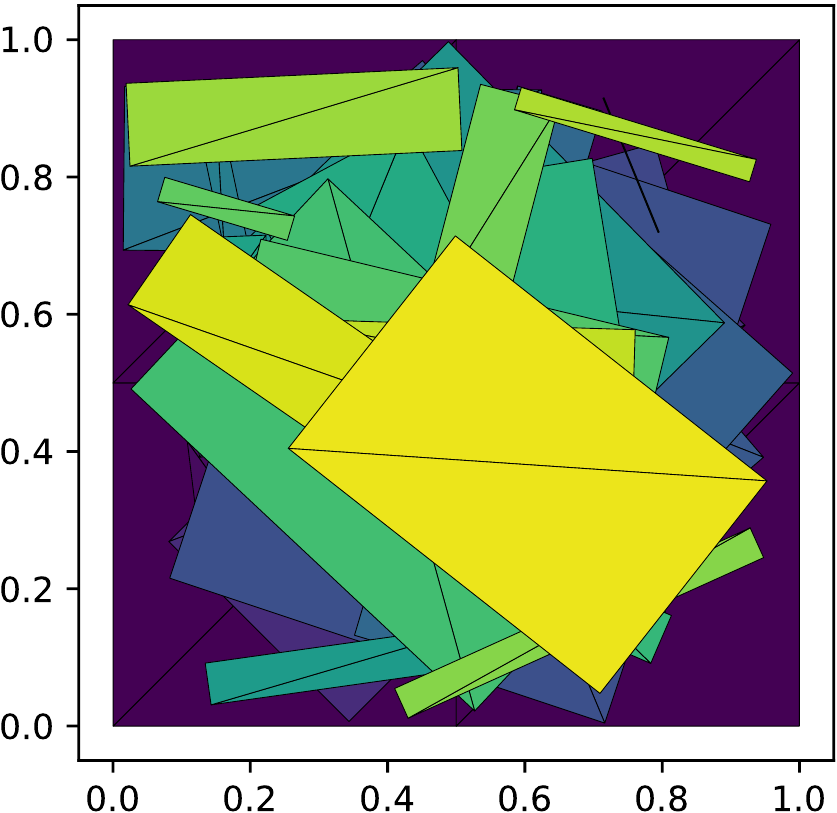}
  \end{center}
  \caption{A sequence of $N$ meshes are randomly placed on top of a fixed background mesh of the unit square for $N=1, 2, 4, 8, 16$ and $32$, shown here for the coarsest refinement level.}
  \label{fig:poisson_meshes}
\end{figure}

Figure~\ref{fig:poisson_meshes} shows the random placement of meshes. Due to the fact that the meshes are randomly placed, some domains may be completely covered by other domains in the hierarchy. For this reason, in the case of $N=8$ domains, there are only $7$ active domains because one of these domains is completely covered. The same holds for $N=16$, where only $15$ of the domains are visible. For the case with $N=32$ domains, three of the domains are completely covered.

\begin{figure}
  \begin{center}
    \includegraphics[width=0.49\textwidth]{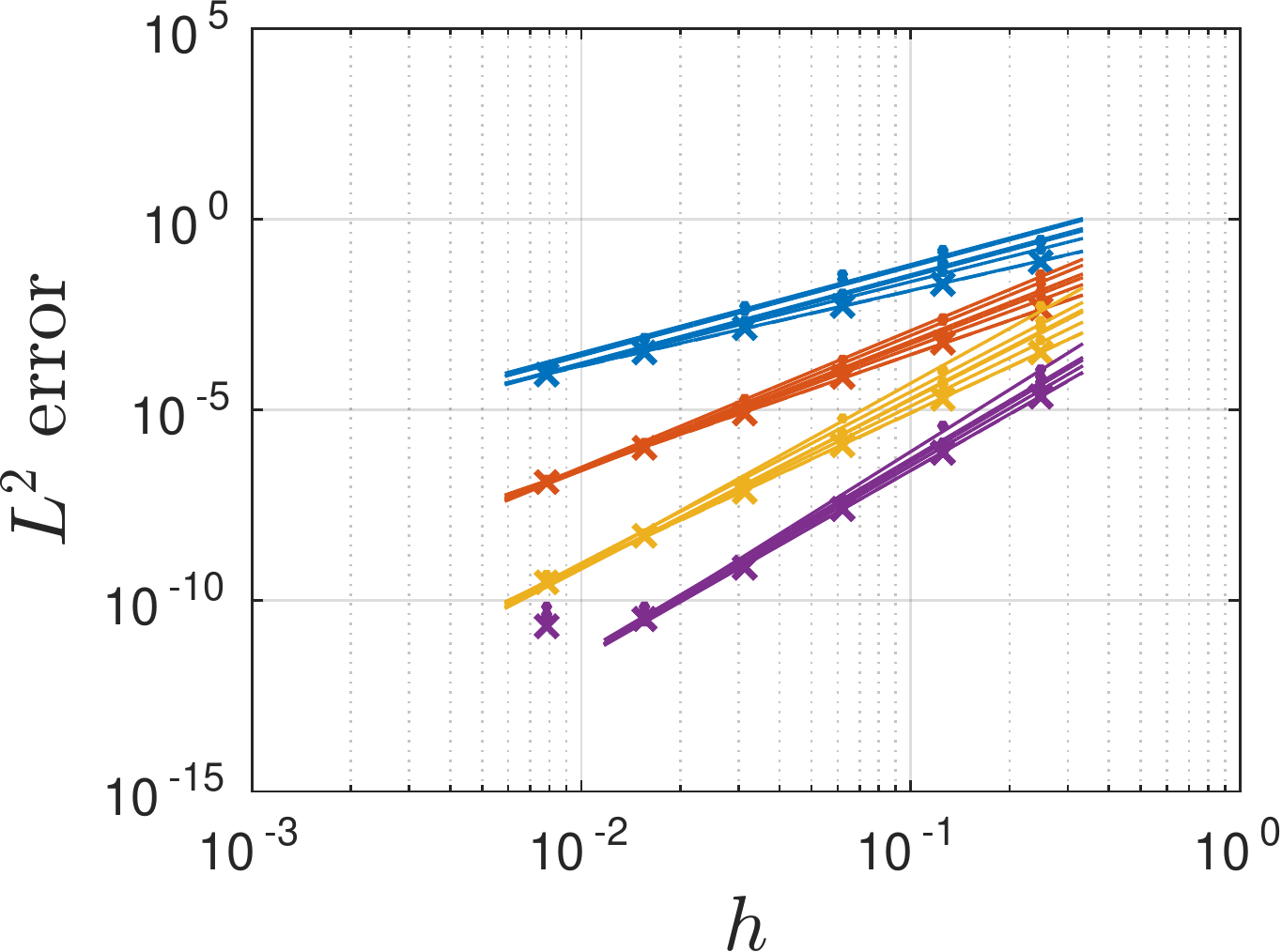}
    \includegraphics[width=0.49\textwidth]{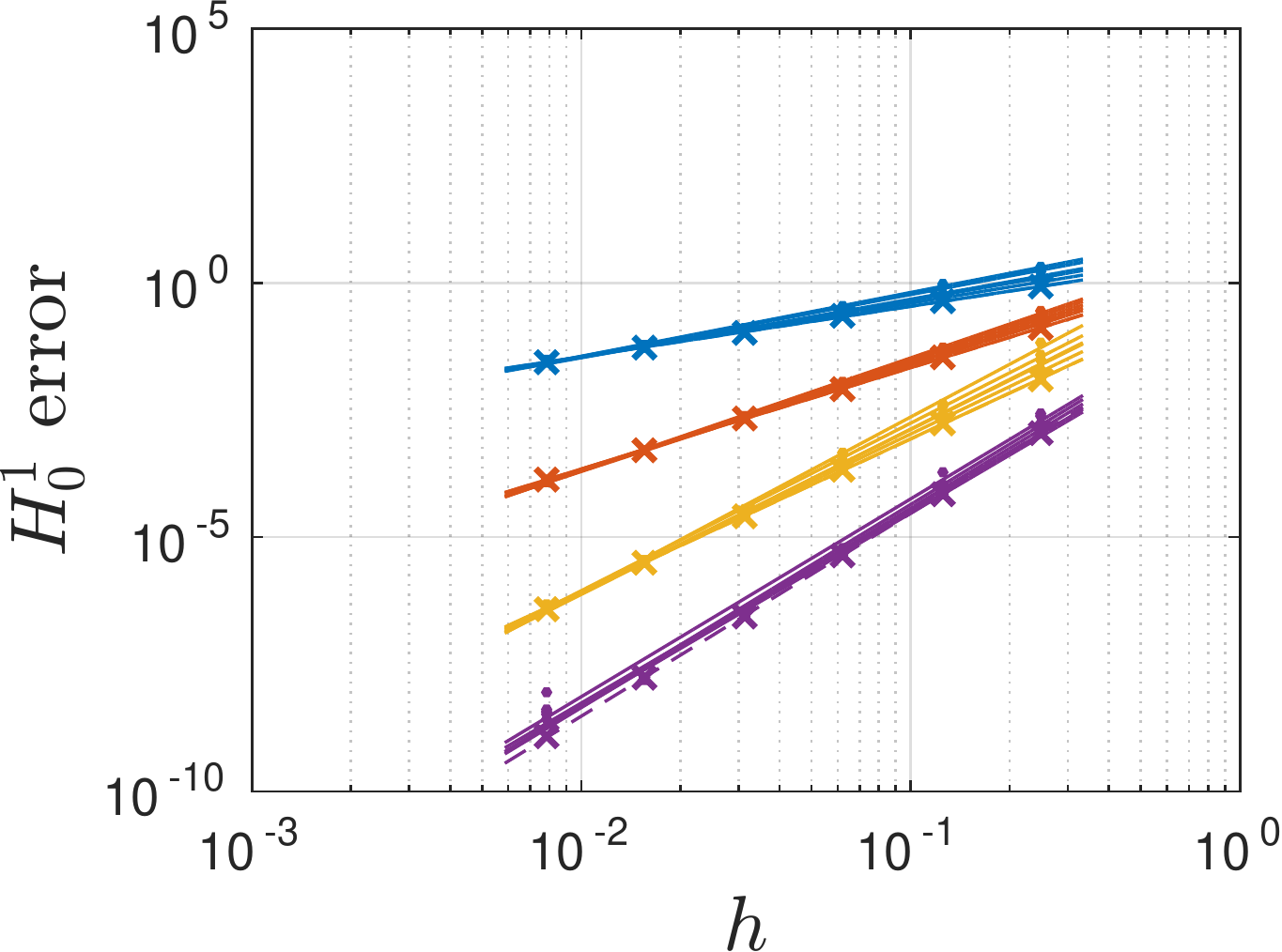}%
  \end{center}
  \caption{Rate of convergence in the $L^2(\Omega)$ (left) and $H^1_0(\Omega)$ (right) norms for $p=1$ (blue), $p=2$ (red), $p=3$
    (yellow) and $p=4$ (purple). For each $p$, the convergence rate is shown for $N = 1, 2, 4, 6, 16, 32$ (six lines) and the errors for $N = 0$ (the standard single mesh discretization) are marked with $\times$ and dashed lines.}
  \label{fig:poisson_convergence}
\end{figure}

Figure~\ref{fig:poisson_convergence} shows the optimal convergence for
polynomial orders $p\in\{1, \ldots, 4\}$ in the $L^2(\Omega)$ and
$H^1_0(\Omega)$ norms. Detailed convergence rates are summarized in Table~\ref{table:poisson_convergence}. As expected, the rate of convergence is $p + 1$ in the $L^2(\Omega)$ norm and $p$ in the $H^1_0(\Omega)$ norm.

We note that the errors are of the same order, independently of the number of meshes. In particular, the error for $N > 0$ subdomains randomly embedded in the unit square is of the same order as for $N = 0$ corresponding to a standard finite element discretization on a single mesh.

A close inspection of Table~\ref{table:poisson_convergence} also reveals that the rate of convergence seems to increase slightly with the number of meshes. The authors believe that this is due to the fact that a larger number of meshes results in a larger number of cut cells and thus an effectively smaller mesh size (more degrees of freedom).

\begin{table}
  \caption{Table with $L^2(\Omega)$ and $H^1_0(\Omega)$ error rates.}
  \label{table:poisson_convergence}
  \centering
  \begin{tabular}{c|c c|c c|c c|c c|}
    \toprule
    &
    \multicolumn{2}{|c|}{$p=1$} &
    \multicolumn{2}{|c|}{$p=2$} &
    \multicolumn{2}{|c|}{$p=3$} &
    \multicolumn{2}{|c|}{$p=4$} \\
    \midrule
    $N$&  $L^2$ & $H^1_0$ &  $L^2$ & $H^1_0$ &  $L^2$ & $H^1_0$ &  $L^2$ & $H^1_0$ \\
    \midrule
0 & 1.9762 & 0.9917 & 2.9965 & 1.9912 & 4.0257 & 3.0044 & 4.9065 & 3.9881 \\
1 & 1.9737 & 0.9911 & 2.9892 & 1.9918 & 4.0235 & 3.0031 & 4.9273 & 3.8062 \\
2 & 2.1712 & 1.0577 & 3.1886 & 2.0496 & 4.2204 & 3.1060 & 5.0149 & 3.7940 \\
4 & 2.3019 & 1.1209 & 3.3259 & 2.1110 & 4.4638 & 3.2509 & 5.1200 & 3.9209 \\
8 & 2.3185 & 1.1522 & 3.3837 & 2.1403 & 4.4083 & 3.2187 & 5.0995 & 3.9484 \\
16 & 2.3347 & 1.2150 & 3.5285 & 2.1880 & 4.4933 & 3.2945 & 5.4149 & 3.8901 \\
32 & 2.3142 & 1.2545 & 3.5842 & 2.2190 & 4.7834 & 3.4492 & 5.4546 & 3.9280 \\
\bottomrule
\end{tabular}
\end{table}

\subsection{Mesh independence}

\begin{figure}
  \begin{center}
    \includegraphics[width=0.3\textwidth]{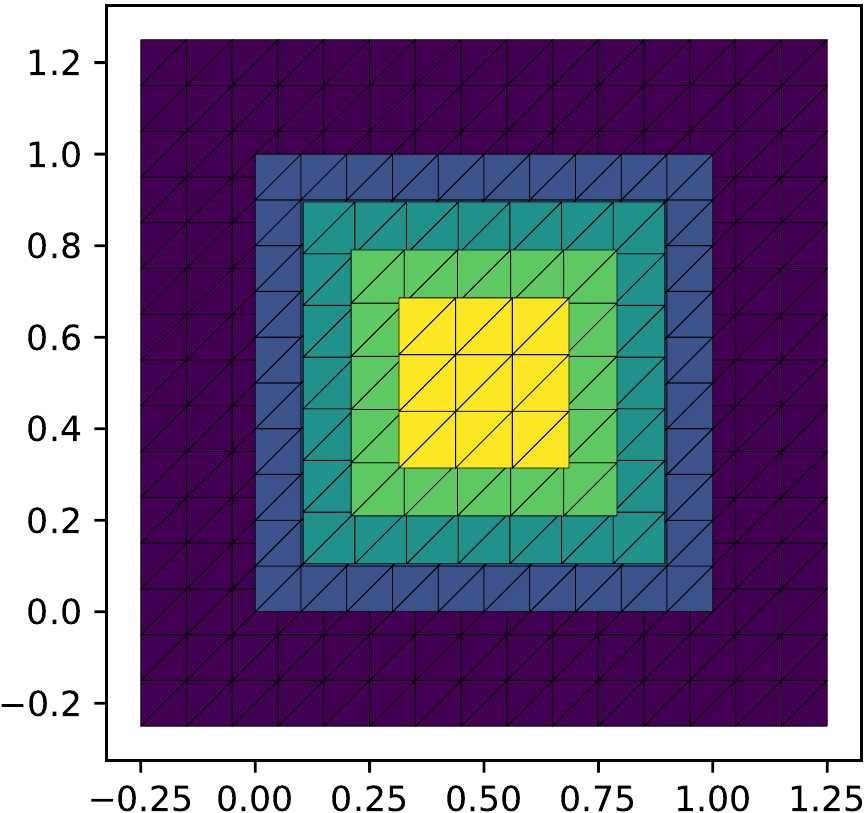}\quad
    \includegraphics[width=0.3\textwidth]{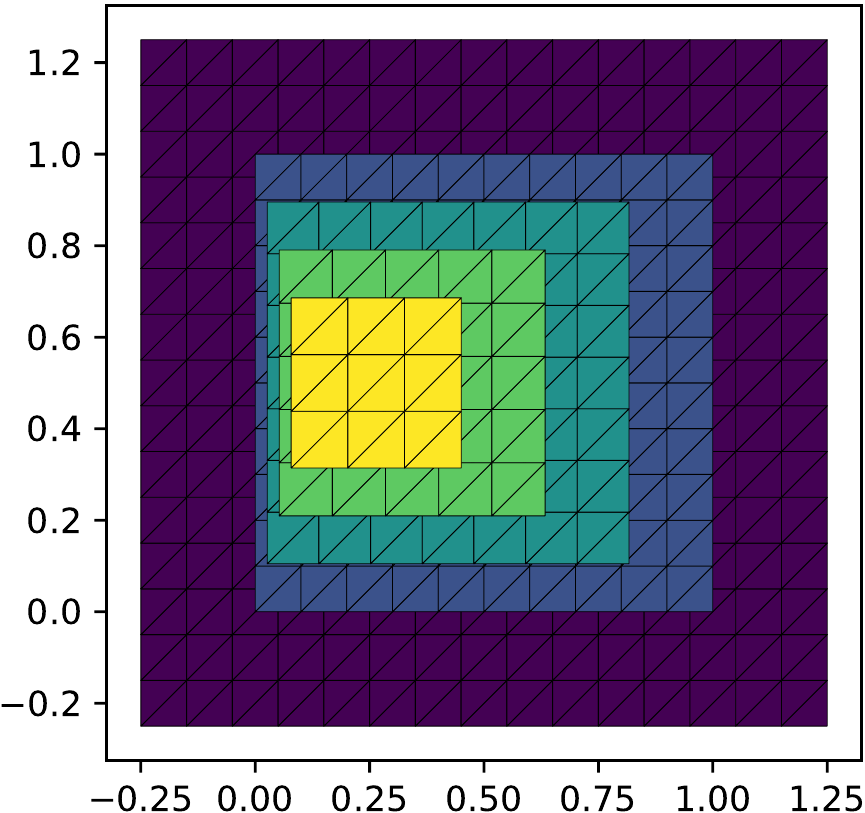}\quad
    \includegraphics[width=0.3\textwidth]{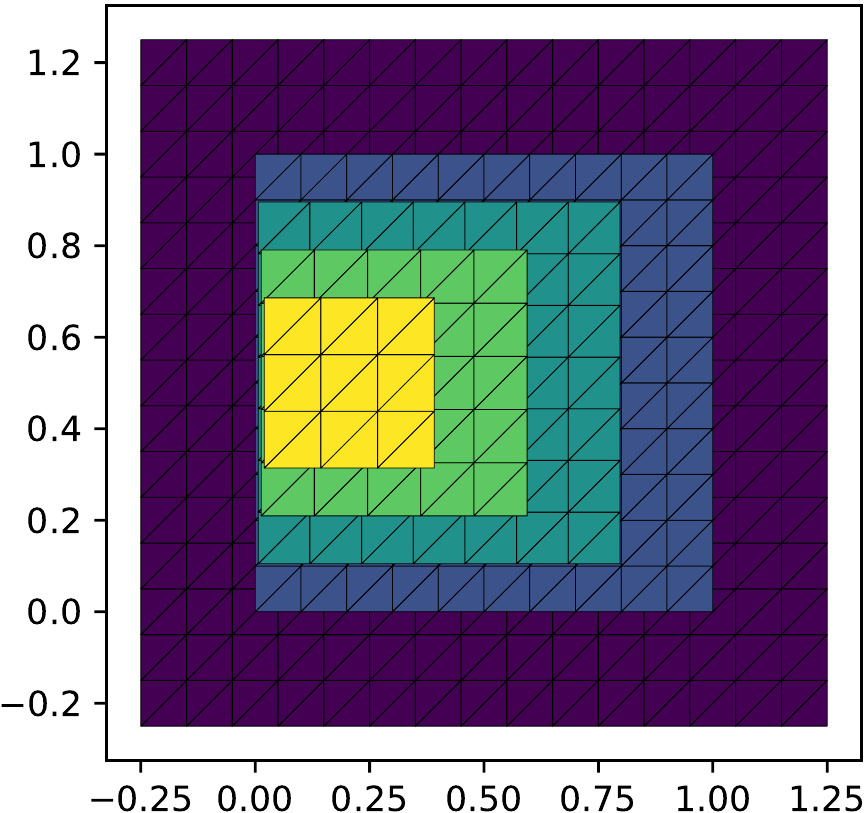}
  \end{center}
  \caption{Mesh configuration for testing the mesh independence of the multimesh formulation in the case of $1 + N = 5$ domains, shown here for
  step $k = 0$, $k = 2$ and $k = 4$.}
 \label{fig:approach_bdry_setup}
\end{figure}

We next investigate the stability of the multimesh finite element method with respect to relative mesh positions, in particular its stability in the presence of thin intersections. We first create a background domain
$\hatOmega_0 = \Omega_0 = [-0.25, 1.25]^2$ and place on top of this the (pre)domain $\hatOmega_1 = [0,1]^2$. We then place $N - 1 = 1, 2, \ldots, 8$
domains inside $\hatOmega_1$ as illustrated in Figure~\ref{fig:approach_bdry_setup} for $N = 4$.

The domains $\{\hatOmega_i\}_{i=2}^N$ are then step-wise moved closer and closer to the left boundary of $\hatOmega_1$ at $x = 0$. In particular, the distance between the left boundaries of the meshes is decreased by a factor $2$ until the distance to $x=0$ is smaller than machine precision $2^{-52} \approx 10^{-16}$. To be precise, the domains $\{\hatOmega_i\}_{i=2}^N$ are defined for $k = 0, 1, \ldots, 52$ as follows:
\begin{align}
  \hatOmega_{i}^k = [x_0^{i,k}, x_0^{i,k}
    + w^i] \times [y_0^i, 1 - y_0^i], \qquad i=2,3,\ldots,N,
\end{align}
where
$a^i = i \pi / (10N)$,
$w^i = 1 - 2a^i$,
$x_0^{i,k} = 2^{-k} a^i$ and
$y_0^i = a^i$.
Note that the shape and size of the domains are kept constant.

For each configuration of the meshes, we compute the multimesh solution of the Poisson problem, using the same exact solution~\eqref{eq:exactsolution} as in the previous example. We use $P_1$ elements with $\beta_0 = 10$ and $\beta_1=5$. The linear system is solved using the conjugate gradient method in PETSc \cite{petsc-user-ref,petsc-efficient} with the hypre BoomerAMG preconditioner \cite{hypre,hypre-web-page}.

\begin{figure}
  \begin{center}
    \includegraphics[height=0.37\textwidth]{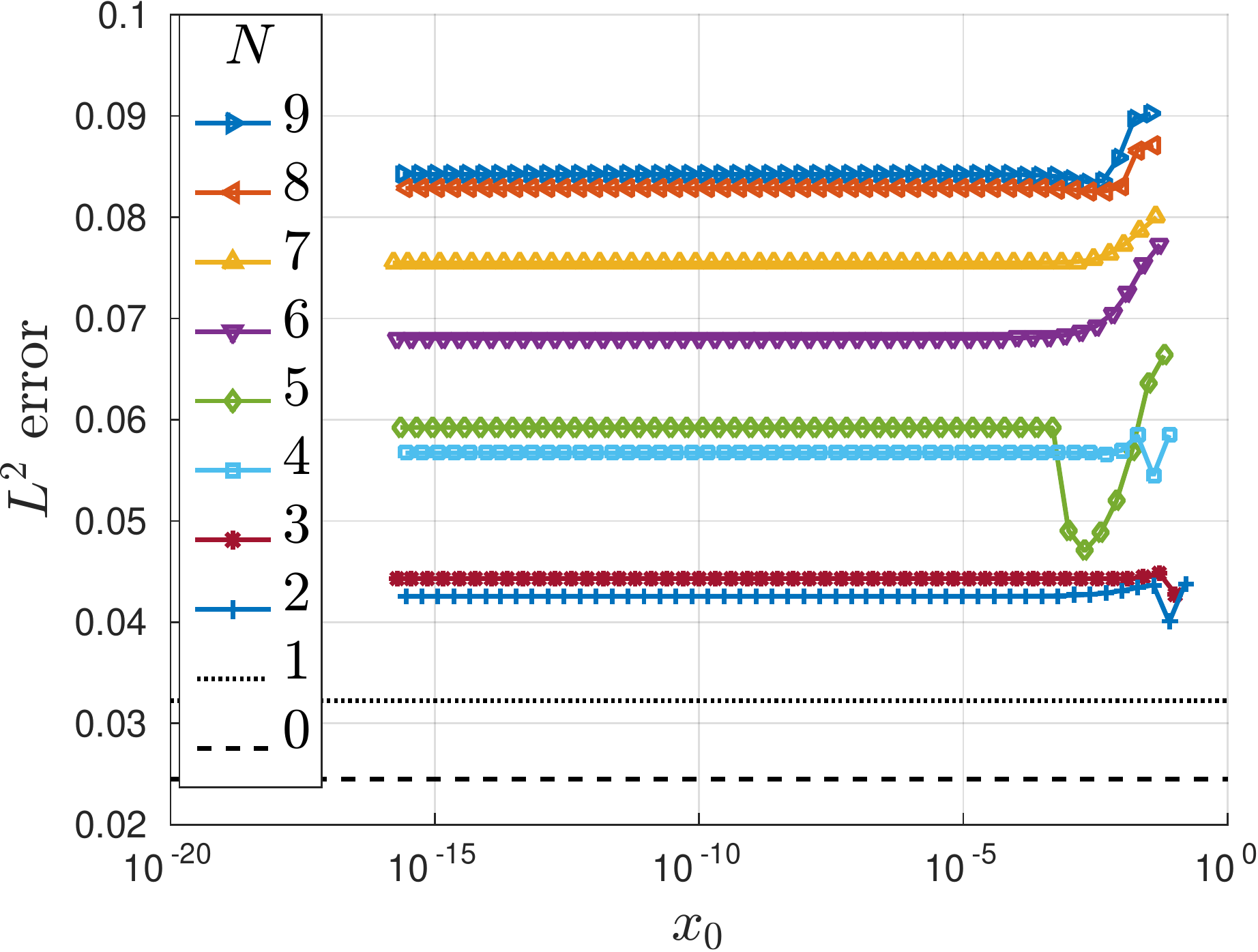}
    \includegraphics[height=0.37\textwidth]{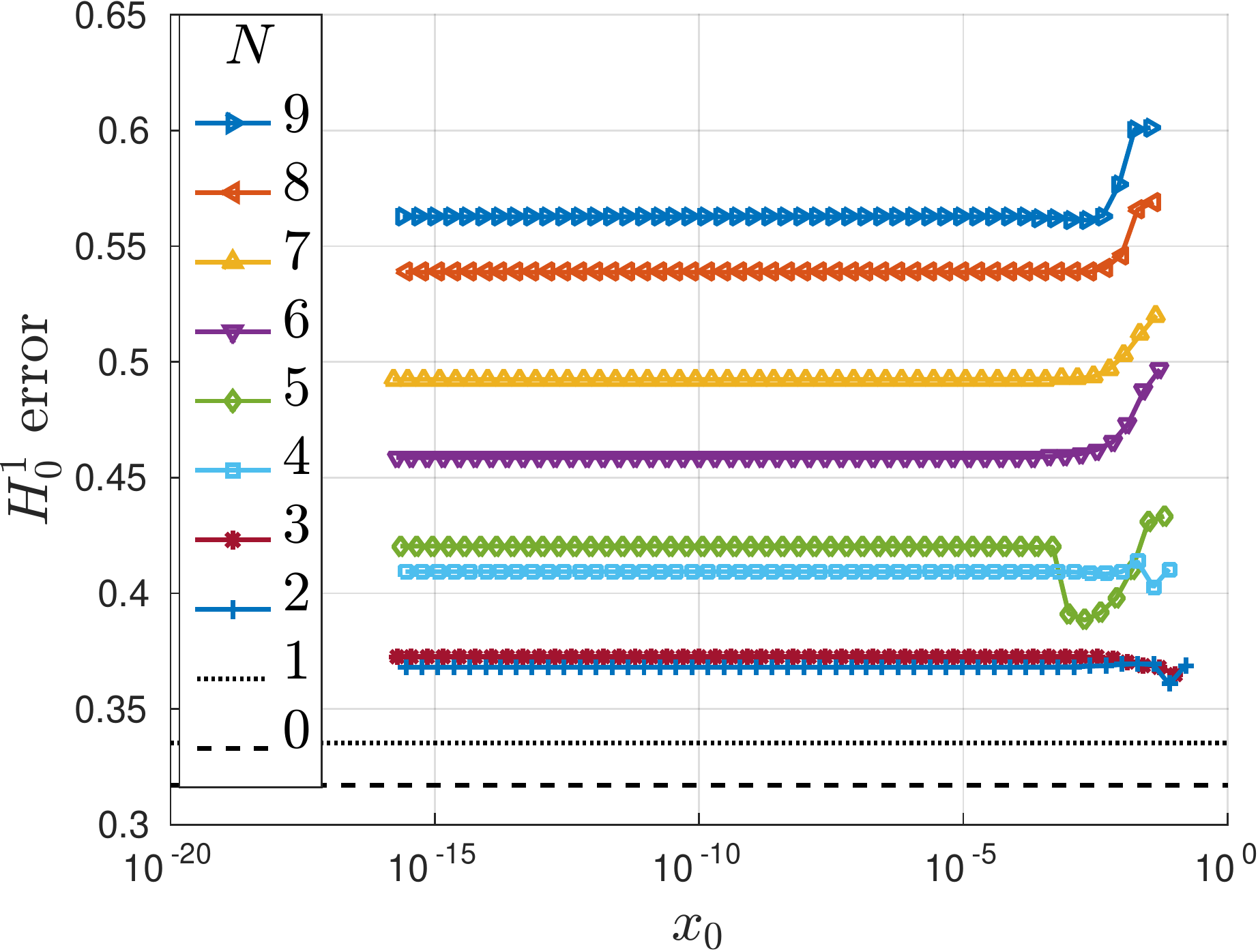}
  \end{center}
  \caption{Errors in the $L^2(\Omega)$ norm (left) and $H^1_0(\Omega)$ norm (right) for varying number of overlapping meshes $N$ as a function of $x_0$. The values for $N=0$ and $N=1$ do not vary with $x_0$ and are only included for comparison.}
  \label{fig:approach_bdry}
\end{figure}

\begin{figure}
  \begin{center}
    \includegraphics[width=0.49\textwidth]{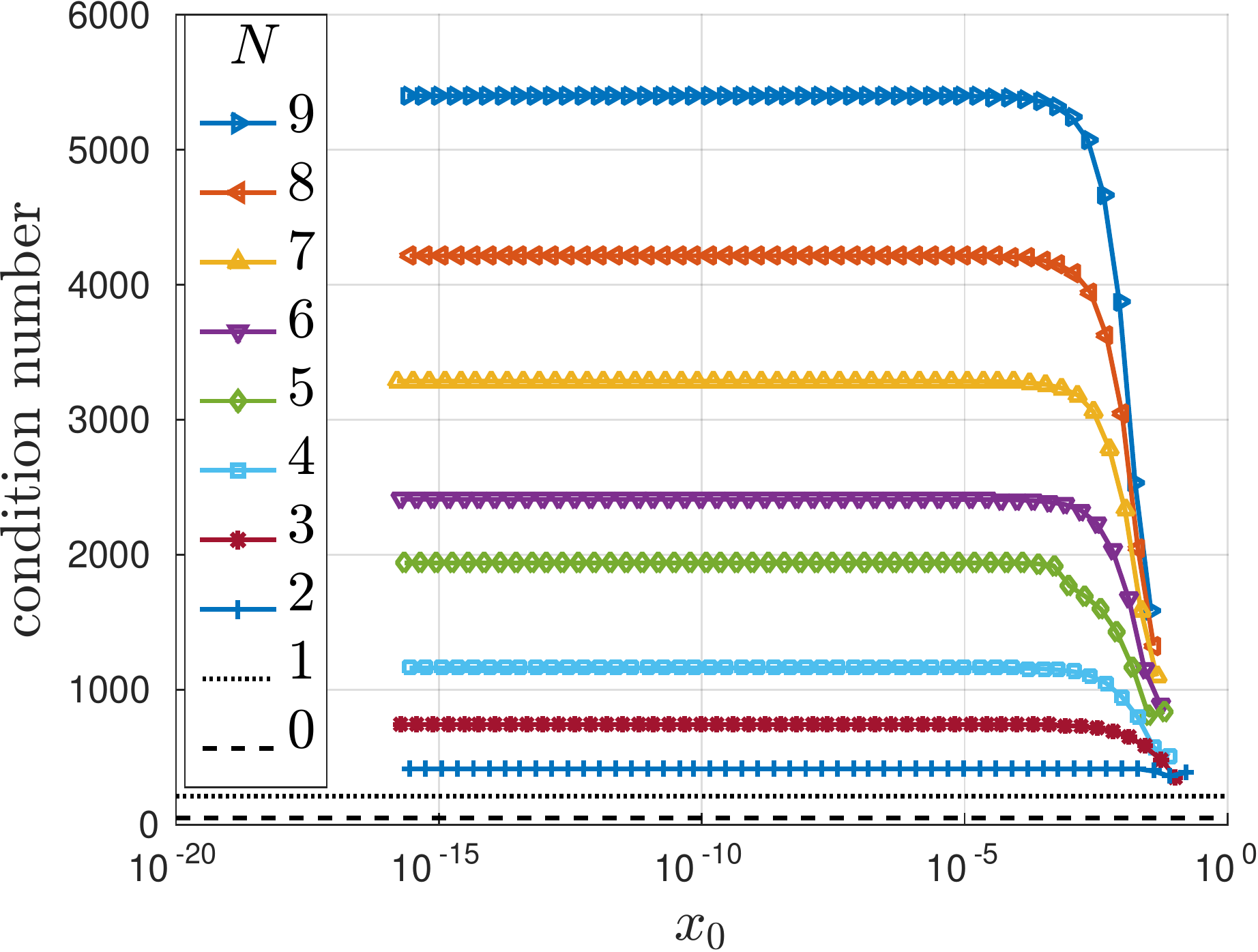}
    \includegraphics[width=0.49\textwidth]{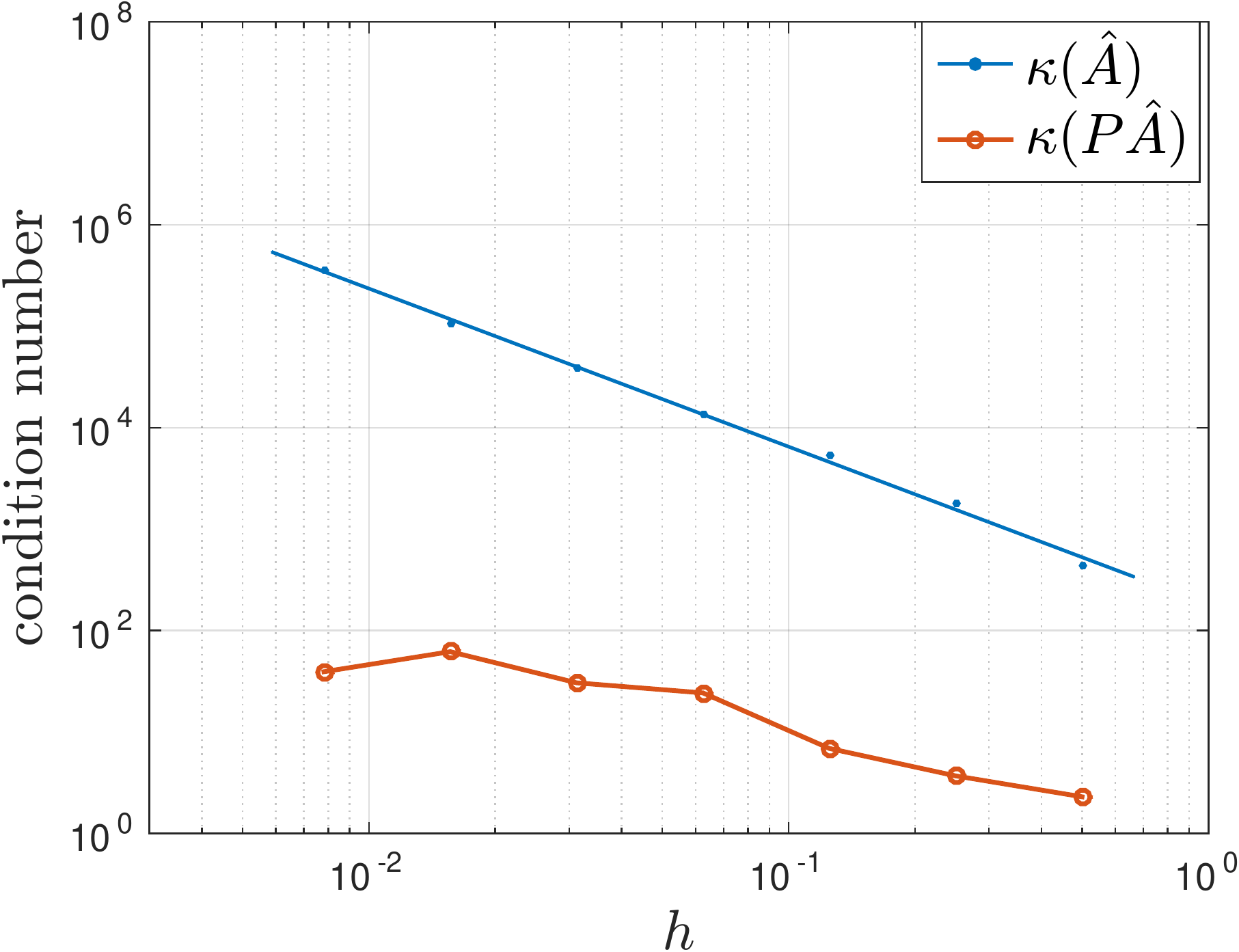}
  \end{center}
  \caption{Left: Condition number as a function of $x_0$ for varying number of overlapping meshes $N$. The values for $N=0$ and $N=1$ do not vary with $x_0$ and are only included for comparison. Right: The condition numbers for $\hat{A}$ and the preconditioned matrix  $P\hat{A}$ as a function of $h$.}
  \label{fig:approach_bdry_cond}
\end{figure}

We first note that both the $L^2(\Omega)$ and $H^1_0(\Omega)$ errors in Figure~\ref{fig:approach_bdry} as well as the condition number in Figure~\ref{fig:approach_bdry_cond} (left) increase slightly with the number of meshes $N$. But most importantly, they all remain constant with increasing $k$, i.e., with increasingly thin mesh intersections. This holds even in the extreme regime when $x_0$ approaches $0$ to within machine precision.

To investigate the condition number we take the ``worst case'', that is, using $N=9$, fix $x_0$ to $x_0 = 2.48\cdot 10^{-16}$ and vary $h$.  $\kappa(\hat{A})$ is computed using the quotient of the largest and smallest singular values as found by SLEPc \cite{slepc-toms,slepc-manual}. The condition number for the preconditioned system matrix, $\kappa(P\hat{A})$, is found from the Krylov solver iterations (the solver is initialized with a random vector and the relative tolerance for convergence is set to $10^{-15}$).

The results are found in Figure~\ref{fig:approach_bdry_cond} (right) and show that $\kappa(\hat{A}) \lesssim  h^{-1.76}$. The last four points have slope $-1.91$ why we believe that $\kappa(\hat{A}) \lesssim  h^{-2}$ in the limit $h\rightarrow 0$. This is the result obtained by the analysis presented in~\cite{mmfem-2} and is the usual result for second-order elliptic problems. The results for $\kappa(P\hat{A})$ seem to be dependent a little on $h$, and the actual values are somewhat large. More work on how to effectively precondition these systems are thus of great interest.

\subsection{Electrostatic interaction of charged bodies}
\label{sec:moving_letters}

As a challenging application of the multimesh formulation and implementation, we model a collection of charged rigid bodies interacting under the influence of the electric field generated by the charges. For simplicity, we set all physical units and material parameters -- masses, mass densities, the electric permittivity, etc.\ -- to $1$.

Consider three disjoint domains $F$, $E$ and $M$ in $\R^2$. The charge density is $\rho = 10$ inside each of the three domains and $0$ outside. In other words,
\begin{equation}
  \rho(x) =
  \begin{cases}
    10, \quad &x \in F \cup E \cup M, \\
    0, \quad  &x \in \R^2 \setminus (F \cup E \cup M).
  \end{cases}
\end{equation}
The electric potential $\phi = \phi(x)$ is then given by the solution of the Poisson problem
\begin{alignat}{2}
  \label{eq:potentialpoisson,1}
  -\Delta \phi &= \rho \quad && \text{in } \R^2, \\
  \label{eq:potentialpoisson,2}
  \phi &= 0 \quad && \text{at } \infty,
\end{alignat}
and the electric field is $E = -\nabla\phi$.

The electric field gives rise to a force $F_D = \int_D E \dx$ and a torque $T_D = \int_D (x - \bar{x}_D) \times E \dx$ on each domain $D = F, E, M$ where $\bar{x}_D$ is the center of mass of the domain $D$. The electric field and the torque generate an acceleration of the center of mass and an angular acceleration around the center of mass for each domain. The three domains are enclosed inside a square domain and contact with the boundary of the box is modeled as a fully elastic collision, i.e., by reversing the normal velocity and angular momentum.

\begin{figure}
  \centering
  \includegraphics[height=0.45\textwidth]{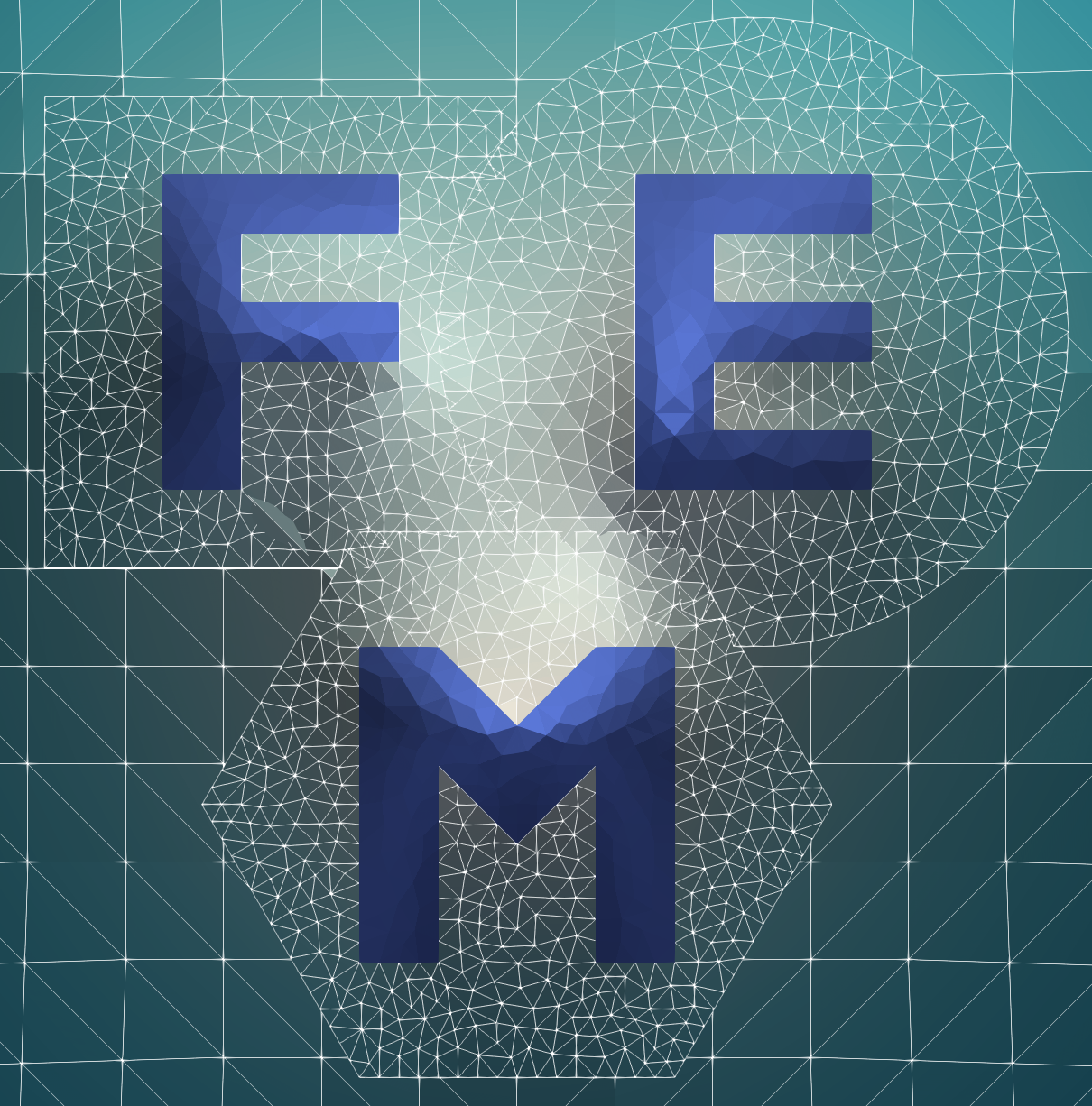} \\[1em]
  \includegraphics[height=0.45\textwidth]{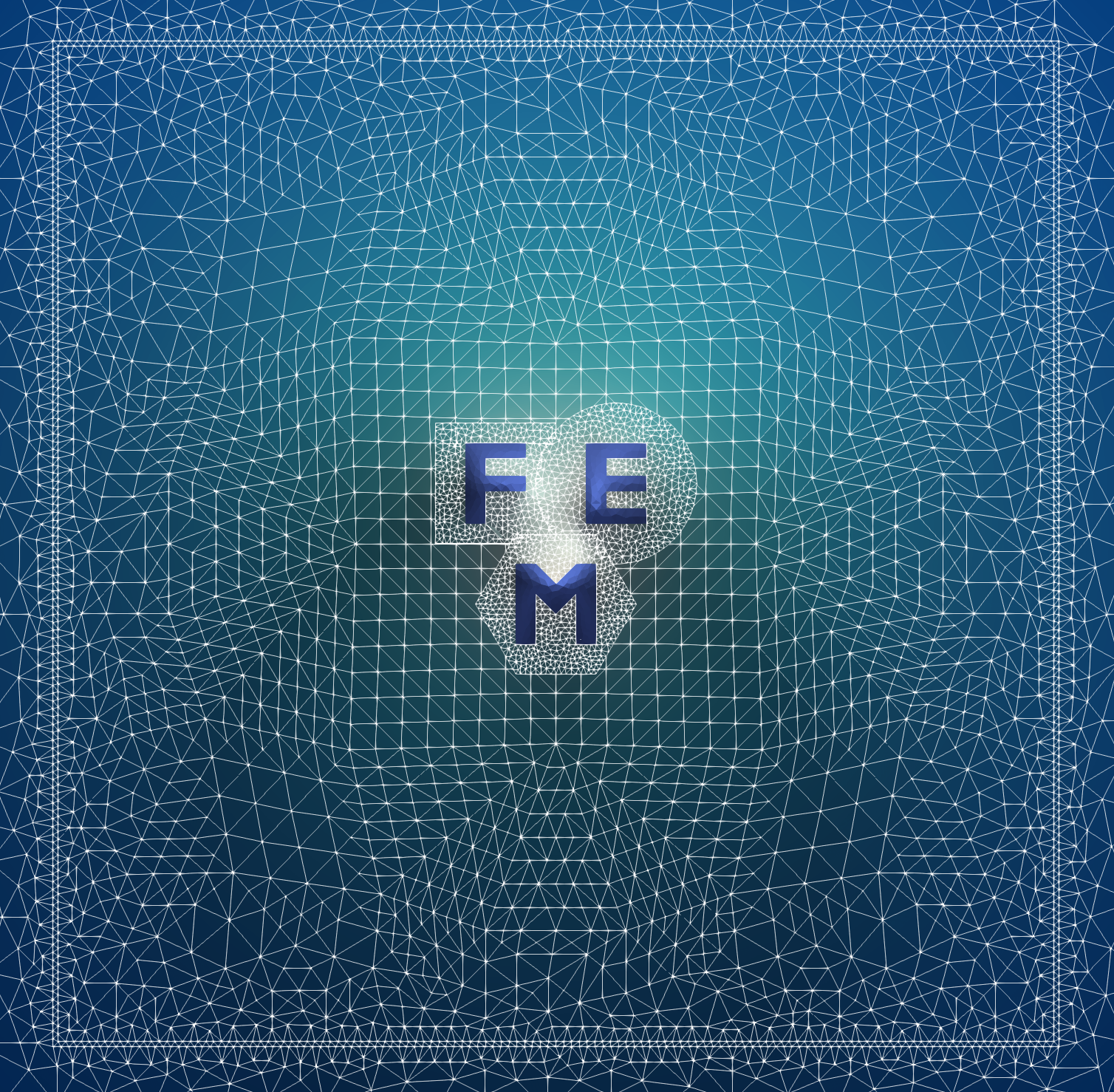} \quad
  \includegraphics[height=0.45\textwidth]{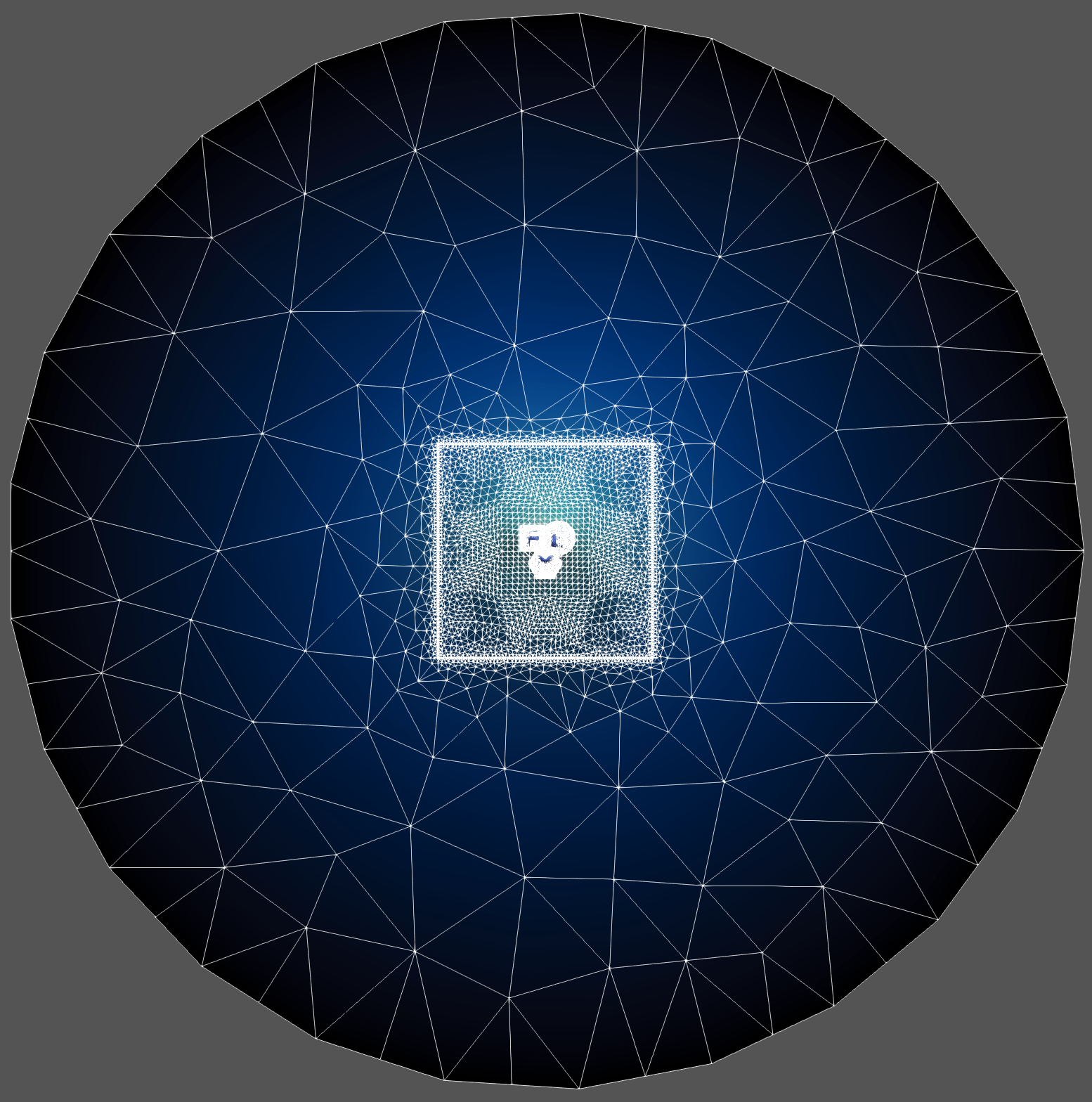}
  \caption{Three domains $F$, $E$ and $M$ modeled as subdomains of three meshes conforming to the domain boundaries. The three meshes are superimposed on a fixed background mesh. The top figure shows a zoom of the three overlapping meshes with $M$ overlapping $E$, and $E$ overlapping $F$ at time $t = 0$. The two bottom figures show two other zooms of the same initial configuration.}
  \label{fig:moving_letters_meshes}
\end{figure}

We simulate the electrostatic interaction of the three bodies by enclosing each domain in a surrounding boundary-fitted mesh and superimposing the three meshes onto a fixed background mesh. The background mesh is made large to model the Dirichlet boundary condition at infinity, and adaptively refined in the interior of the box enclosing the charged bodies, as shown in Figure~\ref{fig:moving_letters_meshes}. We thus have $1 + N = 4$ premeshes, with mesh $\hatmcK_{h,0}$ being the background mesh, and $\hatmcK_{h,i}$, $i=1,2,3$, being the three overlapping meshes. For each overlapping premesh $\hatmcK_{h,i}$, we divide the elements into two subdomains (using markers) to indicate the interior containing the charge and the exterior where the charge is zero.

In each time step, the electric potential $\phi$ is computed using the multimesh discretization of the Poisson problem~\eqref{eq:potentialpoisson,1}--\eqref{eq:potentialpoisson,2}. For sake of variation, we use the alternative stabilization term
\begin{equation}
    s_h (v, w) = \sum_{i=0}^{N-1} \sum_{j=i+1}^N \beta_1 h^{-2} ([ v ], [w ])_{\OO_{ij}}
\end{equation}
with $\beta_1 = 1$; see~\cite{mmfem-2} for a discussion.
The interior penalty parameter is chosen to be $\beta_0 = 10$.

Once the electric potential $\phi_h$ has been computed, the electric field $E_h$ is  computed by projecting the negative gradient of the potential $\phi_h$ to a vector-valued piecewise linear multimesh space. For the projection, we use a similar stabilized variational formulation as for the Poisson problem:
\begin{multline}
   \sum_{i=0}^N (E_h, v_h)_{\Omega_i} +
   \sum_{i=1}^N \sum_{j=0}^{i-1} \beta_0 (h_i^{-1} + h_j^{-1}) ([E_h], [v_h])_{\Gamma_{ij}} \\ +
   \sum_{i=0}^{N-1} \sum_{j=i+1}^N \beta_1 (h_i^{-2} + h_j^{-2})([E_h], [v_h])_{\OO_{ij}}
   =
   \sum_{i=0}^N (-\nabla \phi_h, v_h)_{\Omega_i},
\end{multline}
with $\beta_0 = 10$ and $\beta_1 = 1$. From the electric field, we may then compute the forces and torques on the three rigid bodies $F$, $E$ and $M$.

The dynamics of the rigid bodies is computed using the symplectic Euler scheme; i.e., by first updating the velocities and angular velocities and then updating the positions and rotation angles using the newly computed velocities and angular velocities.

\begin{figure}
  \centering
  \includegraphics[height=0.3\textwidth]{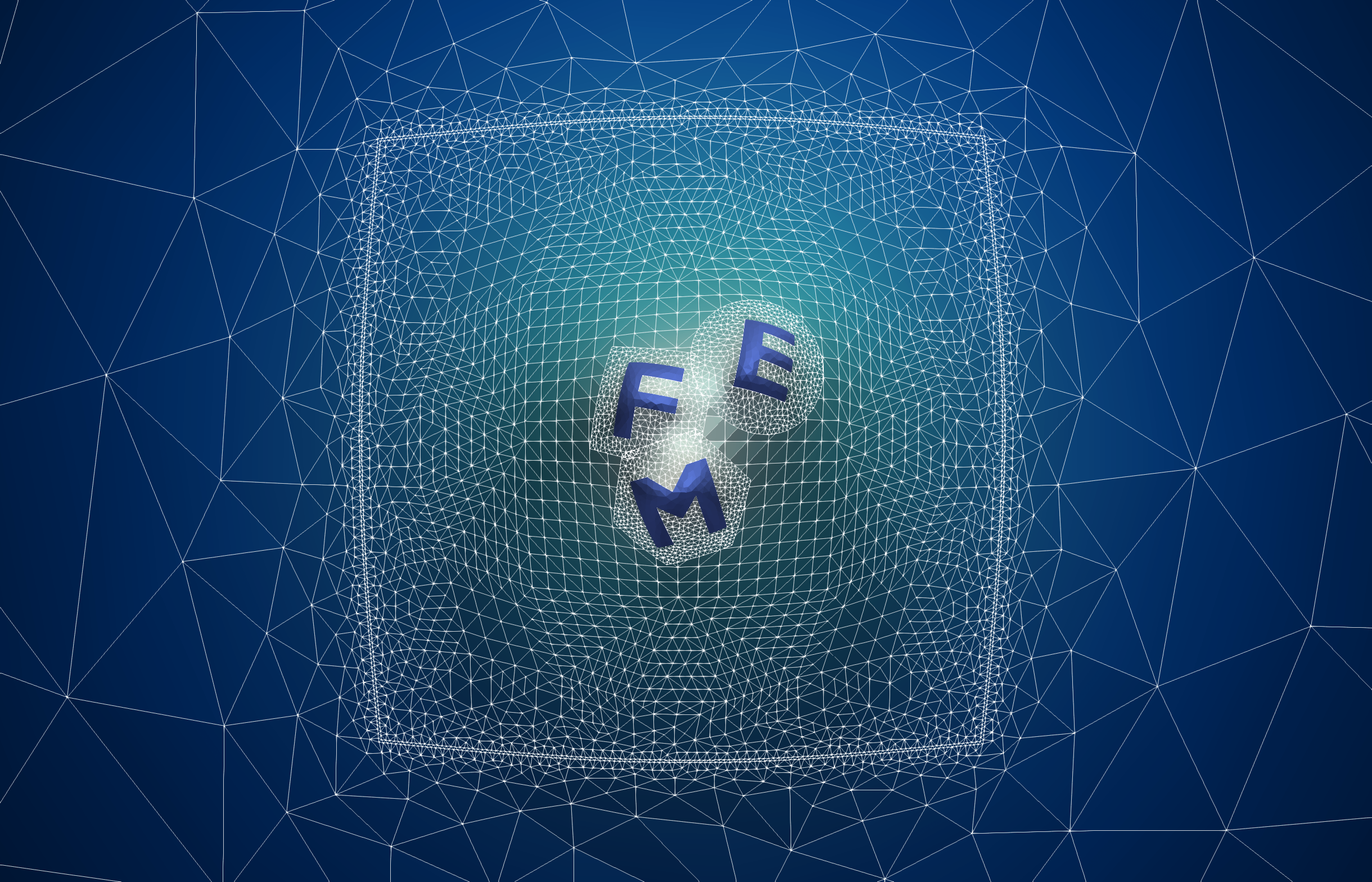} \quad
  \includegraphics[height=0.3\textwidth]{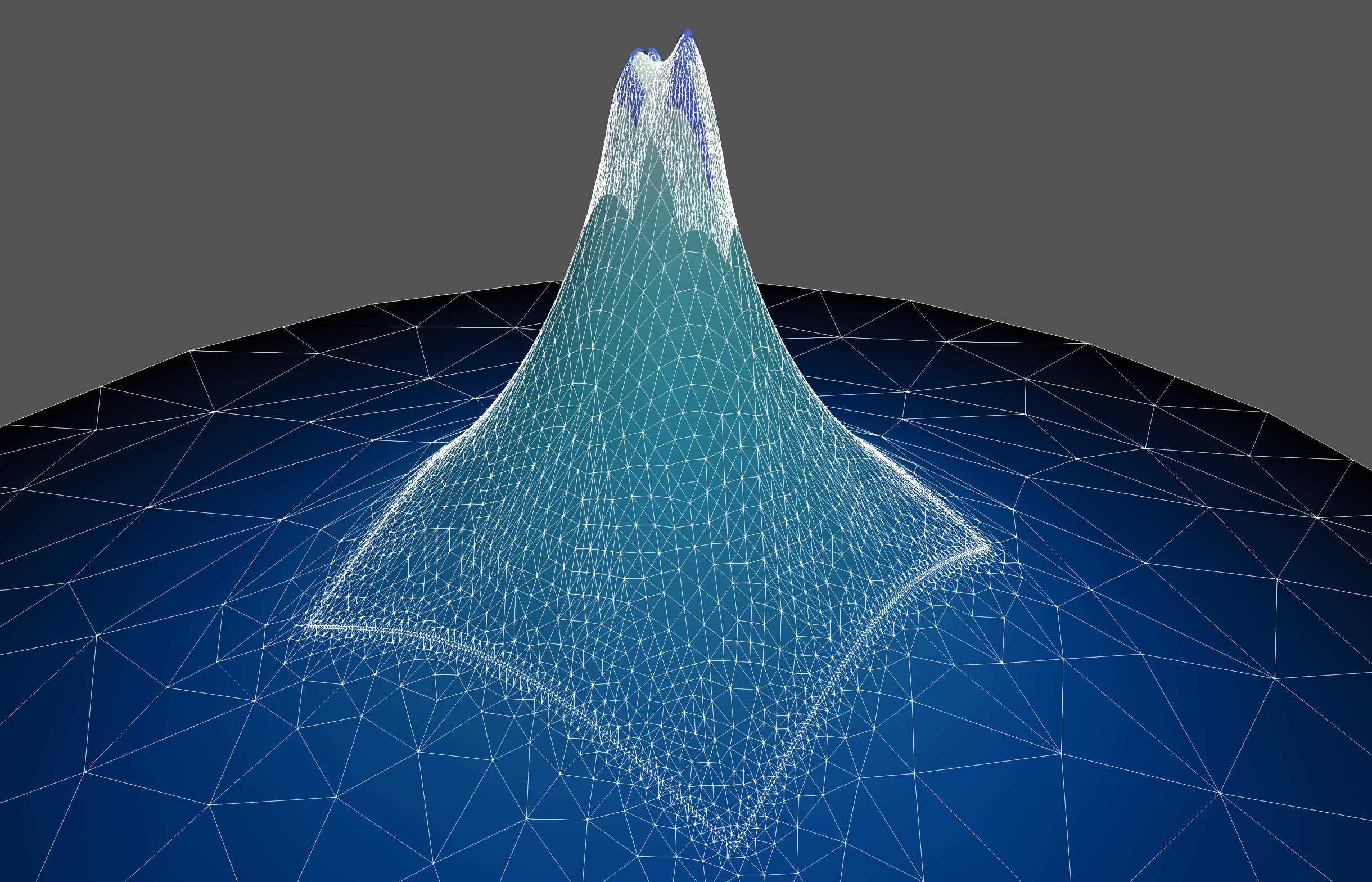}  {}
  \caption{The electric potential $\phi$ at time $t = 9.2$ (time step $46$). The electric field is at each point $x$ described by a finite element representation on the topmost mesh at the point $x$. The finite element representations are glued together via Nitsche terms on the domain boundaries.}
  \label{fig:moving_letters_potential}
\end{figure}

\begin{figure}
  \centering
  \includegraphics[width=0.8\textwidth]{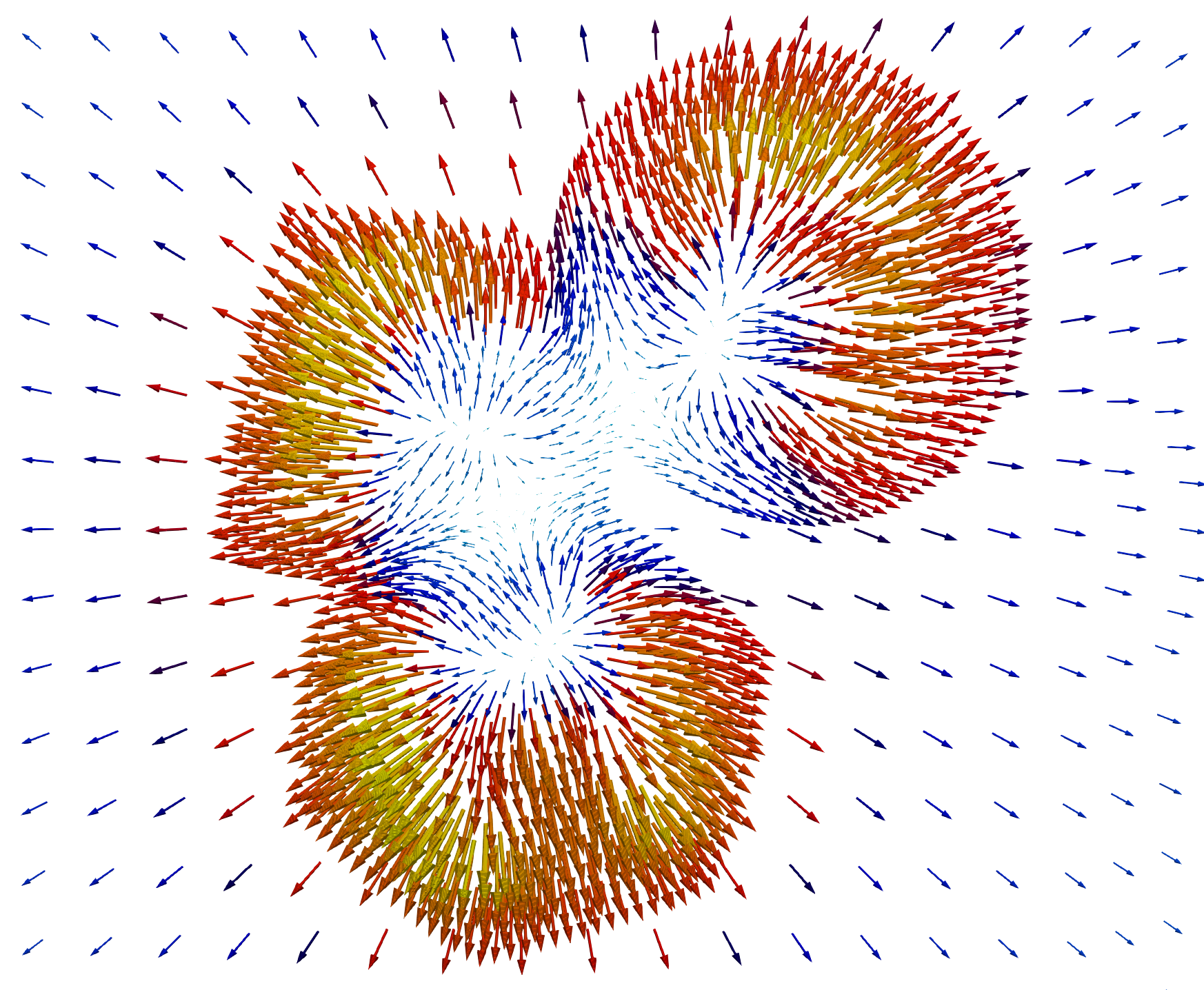}
  \caption{The electric field $E$ at time $t = 9.2$ (time step $46$). Just as the electric potential, the electric field transitions smoothly between the different finite element representations on the three meshes.}
  \label{fig:moving_letters_field}
\end{figure}

The electric charge distribution gives rise to an electric force repelling the three rigid bodies, which then bounce against the boundaries of the enclosing rigid box for 1000 time steps of size $0.2$ ($T = 200$). In each time step, new intersections and quadrature points are computed for the multimesh discretization on the current configuration of the three meshes on top of the fixed background mesh. The computed electric potential $\phi$ and electric field $E$ are shown in Figures~\ref{fig:moving_letters_potential} and~\ref{fig:moving_letters_field}. Throughout the 1000 time steps of the simulation, the electric potential and electric field exhibit a smooth transition across the boundaries of the four meshes, which demonstrates the stability of the finite element formulation and the robustness of the implementation. An animation of the time evolution of the electric potential and the dynamics of the rigid bodies is available on YouTube at \url{https://youtu.be/jBINWmALrps}

%---------------------------------------------------------------------------
\section{Conclusions and Future Work}
\label{sec:conclusions}

We have presented a general framework for discretization of partial differential equations posed on a domain defined by an arbitrary number of intersecting meshes. The framework has been formulated in the context of the Poisson problem, and numerical results presented to support the optimal order convergence and demonstrate the capabilities of the algorithms and implementation. In the accompanying paper~\cite{mmfem-2}, we return to the analysis of the proposed finite element formulation and present optimal order estimates of errors and condition numbers.

In related work~\cite{Johansson:2017ab}, the algorithms and implementation are discussed in more detail. Ongoing and future work by the authors are focused on optimizations and generalizations of the presented formulation, and algorithms and software to enable application to more demanding problems, in particular to multiphysics problems posed on complex and evolving geometries in 3D. The proposed technique for constructing the quadrature rules is dimension independent. However, the collision detection and intersection construction depend on the geometrical (and topological) dimension. The geometric predicate library \cite{Shewchuk:1997aa} provides routines in both 2D and 3D. The intersection construction is currently done in a case-by-case basis and is work in progress. An application to shape optimization is presented in~\cite{Dokken:2017aa}.

%---------------------------------------------------------------------------
\section{Acknowledgements}

August Johansson was supported by The Research Council of Norway through a Centres of Excellence grant to the Center for Biomedical Computing at Simula Research Laboratory, project number 179578, as well as by the Research Council of Norway through the FRIPRO Program at Simula Research Laboratory, project number 25123. Benjamin Kehlet was  supported by The Research Council of Norway through a Centres of Excellence grant to the Center for Biomedical Computing at Simula Research Laboratory, project number 179578.  Mats G.\ Larson was supported in part by the Swedish Foundation for Strategic Research Grant No.\ AM13-0029, the Swedish Research Council Grants Nos.\  2013-4708, 2017-03911, and the Swedish Research Programme Essence. Anders Logg was supported by the Swedish Research Council Grant No.\ 2014-6093.

%---------------------------------------------------------------------------
\bibliographystyle{elsarticle-num}
\bibliography{bibliography}

\end{document}